\documentclass[11pt]{article}

\usepackage{amsmath}
\usepackage{amsthm}
\usepackage{amssymb}
\usepackage{array}
\usepackage[ruled,vlined]{algorithm2e}
\usepackage[nohead,margin=1.1in]{geometry}
\usepackage{color}
\usepackage{verbatim}
\usepackage{booktabs}
\usepackage{threeparttable}
\usepackage{url}

\newtheorem{theorem}{Theorem}[section]

\newtheorem{lemma}{Lemma}[section]
\newtheorem{corollary}{Corollary}[section]
\newtheorem{assumption}{Assumption}[section]

\newtheorem{remark}{Remark}[section]

\usepackage{graphicx,epstopdf}
\usepackage{subfigure}
\usepackage{subfig}
\graphicspath{{./pics/},{./result_v5/}}
\usepackage[toc,page]{appendix}
\allowdisplaybreaks

\numberwithin{equation}{section}
\newcommand{\E}{\mathbb{E}}

\begin{document}
\numberwithin{equation}{section}
\numberwithin{figure}{section}
\title{On the convergence of  stochastic variance reduced gradient for linear inverse problems\thanks{B. Jin is supported by Hong Kong RGC General Research Fund (14306824) and ANR / Hong Kong RGC Joint Research Scheme (A-CUHK402/24) and a start-up fund from The Chinese University of Hong Kong.}}

\author{Bangti Jin\thanks{Department of Mathematics, The Chinese University of Hong Kong, Shatin, N.T., Hong Kong (email: \texttt{b.jin@cuhk.edu.hk}).} \and Zehui Zhou\thanks{Department of Mathematics, The Chinese University of Hong Kong, Shatin, N.T., Hong Kong (email: \texttt{zehuizhou@cuhk.edu.hk}). Author to whom any correspondence should be addressed. }}
\date{}
\maketitle

\begin{abstract}
Stochastic variance reduced gradient (SVRG) is an accelerated version of stochastic gradient descent based on variance reduction, and is promising for solving large-scale inverse problems.
In this work, we analyze SVRG and a regularized version that incorporates \textit{a priori} knowledge of the problem, for solving linear inverse problems in Hilbert spaces. We prove that, with suitable constant step size schedules and regularity conditions, the regularized SVRG can achieve optimal convergence rates in terms of the noise level without any early stopping rules, provided that the truncation level is chosen suitably, and standard SVRG is also optimal for problems with nonsmooth solutions under \textit{a priori} stopping rules. The analysis is based on an explicit error recursion and suitable \textit{a priori} estimates on the inner loop updates 
with respect to the anchor point. Numerical experiments are provided to complement the theoretical analysis.
\end{abstract}

\section{Introduction}

In this work, we consider stochastic iterative methods for solving linear inverse problems in Hilbert spaces:
\begin{equation}\label{eqn:lininv}
A_\dag x=y_\dag,
\end{equation}
where $A_\dag: X \rightarrow Y =Y_1\times \cdots \times Y_n$ denotes the system operator that represents the data formation mechanism and is given by 
\begin{align*}
A_\dag x := (A_{\dag,1} x,\cdots,A_{\dag,n} x)^T, \quad \forall x\in X,
\end{align*}
with bounded linear operators $A_{\dag,i}: X \rightarrow Y_i$ between Hilbert spaces $X$ and $Y_i$ equipped with norms $\|\cdot\|_X$ and $\|\cdot\|_{Y_i}$, respectively, and the superscript $T$ denoting the vector transpose. The element
$x\in X$ denotes the unknown signal of interest and $y_\dag=(y_{\dag,1},\cdots,y_{\dag,n})^T\in Y$ denotes the exact data, i.e., $y_\dag = A_\dag x_\dag$ with $x_\dag$ being the minimum-norm solution relative to the initial guess $x_0$, cf. \eqref{eqn:min-norm} for the definition of $x_\dag$. 
In practice, we only have access to a noisy version $y^\delta$ of the exact data $y_\dag $, given by
\begin{equation*}
y^\delta =(y^\delta_1,\cdots,y^\delta_n)^T=y_\dag +\xi,
\end{equation*}
where $\xi= (\xi_1,\cdots,\xi_n)^T\in Y$ is the noise in the data with a noise level $\delta =\|\xi\|_Y:=\sqrt{\sum_{i=1}^n\|\xi_i\|_{Y_i}^2}.$ Below we assume $\delta<1$. 
Linear inverse problems arise in many practical applications, e.g., computed tomography \cite{HermanLentLutz:1978,StrohmerVershynin:2009,GaoBlumensath:2018}, {photoacoustic tomography \cite{PoudelLouAnastasio:2019,XiaWangHan:2022}}  and positron emission tomography \cite{HudsonLarkin:1994,KeretaTwyman:2021,TwymanThielemans:2023}.

Stochastic iterative algorithms, including stochastic gradient descent (SGD) \cite{RobbinsMonro:1951,JinLu:2019,JahnJin:2020,JinZhouZou:2021siuq,ZhouDSGD:2024} and stochastic variance reduced gradient (SVRG) \cite{JohnsonZhang:2013,ZhangMahdaviJin:2013,JinZhouZou:2022ip}, have gained much interest in the inverse problems community in recent years, due to their excellent scalability with respect to data size.  We refer interested readers to the recent surveys \cite{EhrhardtKereta:2025,JinXiaZhou:2025} for detailed discussions. 
Specifically, consider the following ill-posed optimization problem
\begin{equation*}
 J (x) = \frac{1}{2n}\|A_\dag x-y^\delta\|_Y^2 = \frac{1}{n}\sum_{i=1}^n f_i(x),\quad \mbox{with}\quad f_i(x)=\frac{1}{2} \|A_{\dag,i} x-y_i^\delta\|_{Y_i}^2.
\end{equation*}
Given an initial guess $x_0^\delta\equiv x_0\in X$, the SGD iteration is given by
\begin{equation*}
x_{k+1}^\delta = x_k^\delta - \eta_k f'_{i_k}({ x_k^\delta}), \quad k=0,1,\cdots,
\end{equation*}
while the SVRG iteration reads
\begin{equation*}
x_{k+1}^\delta =x_k^\delta -\eta_k \big( f'_{i_k}(x_k^\delta)- f'_{i_k}(x_{[\frac{k}{M}]M}^\delta)+ J'(x_{[\frac{k}{M}]M}^\delta)\big), \quad k=0,1,\cdots,
\end{equation*}
where $\{\eta_k\}_{k\geq 0}\subset (0,\infty)$ is the step size schedule, the index $i_k$ is sampled uniformly at random from the index set $\{1,\ldots,n\}$, $M$ is the frequency of computing the full gradient, and the notation $[\cdot]$ denotes taking the integral part of a real number. 

By combining the full gradient $J'(x_{[\frac{k}{M}]M}^\delta)$ of the objective $J$ at the anchor point $x_{[\frac{k}{M}]M}^\delta$ with a random gradient gap $f_{i_k}'(x_k^\delta)-f_{i_k}'(x_{[\frac{k}{M}]M}^\delta)$, 
SVRG can accelerate the convergence of SGD and has become very popular in stochastic optimization \cite{BottouCurtisNocedal:2018,GowerSchmidtBach:2020}. Its performance depends on the frequency $M$ of computing the full gradient, and $M$ was suggested to be $2n$ and $5n$ for convex and nonconvex optimization, respectively \cite{JohnsonZhang:2013}.  In practice,
there are several variants of SVRG, depending on the choice of the anchor point, e.g., last iterate and randomly selected iterate within the inner loop.
In this work, we focus on the version given in Algorithm \ref{alg:svrg}, where $A_\dag^*$ and $A_{\dag,i}^*$ denote the adjoints of the operators $A_\dag$ and $A_{\dag,i}$, respectively.

\medskip
\begin{algorithm}[H]
\SetAlgoLined
Set initial guess $x_0^\delta=x_0$, frequency $M$, and step size schedule $\{\eta_k\}_{k\geq0}$\\
 \For{$K=0,1,\cdots$}{
compute $ {g_K = J'(x_{KM}^\delta)=\tfrac{1}{n}A_\dag^*(A_\dag x_{KM}^\delta-y^\delta) }$\\
\For{$t=0,1,\cdots,M-1$}{
draw {$i_{KM+t}$} i.i.d. uniformly from $\{1,\cdots,n\}$\\
update $x_{KM+t+1}^\delta =x_{KM+t}^\delta -\eta_{KM+t} \big(A_{\dag,i_{KM+t}}^*A_{\dag,i_{KM+t}}(x_{KM+t}^\delta-x_{KM}^\delta)+ g_K\big)$\\
}
check the stopping criterion.
}
\caption{SVRG for problem \eqref{eqn:lininv}.\label{alg:svrg}}
\end{algorithm}

The low-rank nature of  $A_\dag$ implies that one can extract a low-rank subspace. Several works have proposed subspace / low-rank versions of stochastic algorithms \cite{Kozak:2021,Gressmann2020,LiWangFan:2022,LiangLiu::2024,HeLiHu2025}. Let $A:=(A_1, \cdots, A_n)^T$ approximate $A_\dag$. Using $A$ in place of $A_\dag$ in Algorithm \ref{alg:svrg} gives Algorithm \ref{alg:dsvrg}, termed as regularized SVRG  (rSVRG) below. The version described in Algorithm \ref{alg:dsvrg} follows Algorithm \ref{alg:svrg} directly, and there are also other variants. One such variant based on an approximate singular value decomposition (SVD) $A(\cdot) = \sum_{j=1}^J \sigma_j (\varphi_j,\cdot)\psi_j$ of $A_\dag$ may take the form $ A_j(\cdot) = \sigma_j(\varphi_j,\cdot)\psi_j$, $j=1,\ldots,J$. This alternative version is denoted by rSVRG$_*$ below. Numerically computing SVD directly can be expensive, although randomized SVD \cite{Halko:2011} may allow extracting low-rank approximations efficiently. See Section \ref{sec:num} for a comparison between rSVRG and rSVRG$_*$.
Conceptually speaking, rSVRG may be interpreted as integrating a learned prior into SVRG, if the approximation $A$ is generated from a paired training dataset $\{x^{(j)},y^{(j)}\}_{j=1}^N$. The regularization provided by the learned prior may relieve the need of early stopping. 

\begin{algorithm}[hbt!]
\SetAlgoLined
Set initial guess $x_0^\delta=x_0$, frequency $M$, and step size schedule $\{\eta_k\}_{k\ge0}$\\
 \For{$K=0,1,\cdots$}{
compute $ g_K = \tfrac{1}{n}A^*(A x_{KM}^\delta-y^\delta)$\\
\For{$t=0,1,\cdots,M-1$}{
draw {$i_{KM+t}$} i.i.d. uniformly from $\{1,\cdots,n\}$\\
update $x_{KM+t+1}^\delta =x_{KM+t}^\delta -\eta_{KM+t} \big(A_{i_{KM+t}}^*A_{i_{KM+t}}(x_{KM+t}^\delta-x_{KM}^\delta)+ g_K\big)$\\
}
check the stopping criterion.
}
\caption{Regularized SVRG (rSVRG)  for problem \eqref{eqn:lininv}.\label{alg:dsvrg}}
\end{algorithm}

The mathematical theory of SVRG for inverse problems from the perspective of regularization theory has not been fully explored, and only recently has its convergence rate for solving linear inverse problems been investigated  \cite{JinZhouZou:2022ip,JinChen:2024}. 
In this work, we establish the convergence of both rSVRG and SVRG for solving linear inverse problems. See Theorem \ref{thm:main} for the error bounds in terms of the iteration index $k$, Corollary \ref{cor:rate} for the convergence rates in terms of $\delta$, and Corollary \ref{cor:regularizing} for regularizing property. Note that rSVRG inherits the regularization mechanism of the truncated SVD, without any need of early stopping rules, provided that the truncation level is chosen properly according to the noise level, and empirically, it can outperform SVRG (i.e., higher accuracy), cf. the numerical results in Section \ref{sec:num}. 
Moreover, we establish the (optimal) convergence rates in both expectation and uniform sense for both SVRG (when combined with \textit{a priori} stopping rules) and rSVRG (cf. Theorem \ref{thm:main} and Corollary \ref{cor:rate}), while the prior works \cite{JinZhouZou:2022ip,JinChen:2024} only studied convergence rates in expectation. Note that the work \cite{JinZhouZou:2022ip} is based on the spectral analysis of bounded linear operators, whereas the work \cite{JinChen:2024} is based on a Lyapunov type analysis. In this work, we follow the framework of the spectral analysis in \cite{JinZhouZou:2022ip}.
For SVRG, the condition for its optimal convergence rate in expectation is more relaxed than that in \cite{JinZhouZou:2022ip}. 
However, unlike the results in \cite{JinZhouZou:2022ip}, SVRG loses its optimality for smooth solutions under the relaxed condition.
For the benchmark source condition $x_\dag-x_0\in\mathrm{Range}(A_\dag^*)$ studied in \cite{JinChen:2024}, the  condition is either comparable  or more relaxed; see Remark \ref{rem:rate} for further details.

The rest of the work is organized as follows. In Section \ref{sec:main}, we present and discuss the main results, i.e., the convergence rate for (r)SVRG in Theorem \ref{thm:main} and the regularizing properties in Corollary \ref{cor:regularizing}. We present the detailed proofs in Section \ref{sec:conv}. Then in Section \ref{sec:num}, we present 
several numerical experiments to complement the analysis. 
Finally, we conclude this work with further discussions in Section \ref{sec:conc}. 
In Appendix \ref{app:estimate}, we collect lengthy and technical proofs of several technical results. Throughout, we suppress the subscripts in the norms and inner products, as the spaces are clear from the context.

\section{Main result and discussions}\label{sec:main}
To present the main result of the work, we first state the assumptions on the step size schedule $\{\eta_j\}_{j\geq 0}$, the reference solution $x_\dag$, the unique minimum-norm solution relative to $x_0$, given by
\begin{equation}\label{eqn:min-norm}
  x_\dag = \arg\min_{x\in X: A_\dag x = y_\dag} \|x-x_0\|,
\end{equation}
and the operator $A$, for analyzing the convergence of  rSVRG. We denote the operator norm of $A_i$ by $\|A_i\|$ and that of $A$ by $\|A\|\leq\sqrt{\sum_{i=1}^n\|A_i\|^2}$. $\mathcal{N}(A_\dag)$ denotes the null space of $A_\dag$.

\begin{assumption}\label{ass}
The following assumptions hold.
\begin{itemize}
\item[$\rm(i)$] The step size $\eta_j = c_0$, $j=0,1,\cdots$, with $c_0\leq L^{-1}$, where $L:= \max_{1\leq i\leq n}\|A_i\|^2$.
\item[$\rm(ii)$] There exist $\nu>0$ and $w\in \mathcal{N}(A_\dag)^\perp$ such that $x_\dag-x_0=B_\dag^\nu w$, with $B_\dag=n^{-1}A_\dag^* A_\dag$
and $\mathcal{N}(A_\dag)^\perp$ being the orthogonal complement of $\mathcal{N}(A_\dag)$.
\item[$\rm(iii)$] Let $a\geq0$ be a constant. 
When $a=0$, set $A=A_\dag$. 
When $a>0$, let $A_\dag$ be a compact operator with $\{\sigma_{j}, \varphi_{j}, \psi_{j}\}_{j=1}^\infty$ being its singular values and functions, i.e., $A_\dag(\cdot)=\sum_{j=1}^\infty \sigma_{j}\langle\varphi_{j},\cdot\rangle\psi_{j}$, such that $\{\sigma_{j}\}_{j=1}^{\infty}\subset [0,\infty)$, $\sigma_{j}\geq\sigma_{j'}\geq a\delta^b> 0$ for any $j\leq j'\leq J$, and $\sigma_{j}<a\delta^b$ for any $j>J$, with some $b>0$.
Set $A(\cdot)=\sum_{j=1}^{J} \sigma_{j}\langle\varphi_{j},\cdot\rangle\psi_{j}$.
\end{itemize}
\end{assumption}

The constant step size in Assumption \ref{ass}(i) is commonly employed by SVRG \cite{JohnsonZhang:2013}.  (ii) is commonly known as the source condition \cite{EnglHankeNeubauer:1996}, which imposes certain regularity on the initial error $x_\dag-x_0$ and is crucial for deriving convergence rates for (iterative) regularization methods.
Without the condition, the convergence of regularization methods can be arbitrarily slow \cite{EnglHankeNeubauer:1996}. 
(iii) assumes that the operator $A$ captures the important features of $A_\dag$, and can be obtained by the truncated SVD of $A_\dag$ that retains principal singular values $\sigma_{j}$ such that $\sigma_{j}\geq a\delta^{b}$. 
When $a=0$, $A=A_\dag$ and rSVRG reduces to the standard SVRG. (iii) is adopted to simplify the error analysis, and one can obtain similar results even if it is relaxed, which can be fulfilled by several data-driven models, including data-driven reduced order models \cite{XieMohebujjaman2018,MouKocSan2021}, autoencoder neural networks \cite{Kashima2016}, and neural networks combined with model reduction \cite{BhattacharyaHosseiniStuart:2021,JinZhouZou:2024}. See Corollary \ref{cor:svd0} for details of the relevant convergence results in the absence of (iii).

Let $\mathcal{F}_k$ denote the filtration generated by the random indices $\{i_0,i_1,\ldots,i_{k-1}\}$, $\mathcal{F}={\bigvee_{k=1}^\infty} \mathcal{F}_k$, $(\Omega,\mathcal{F},\mathbb{P})$ denote the associated probability space, and $\E[\cdot]$ denote taking the expectation with respect to the filtration $\mathcal{F}$. 
The (r)SVRG iterate $x_k^\delta$ is random but measurable with respect to the filtration $\mathcal{F}_k$. 
Now, we state the main result on the error $e_k^\delta=x_k^\delta-x_\dag$ of the (r)SVRG iterate $x_k^\delta$ with respect to the minimum-norm solution $x_\dag$. Below we adopt the convention $k^0:=\ln k$, and let 
\begin{align*}
    \overline{C_0}&:=\max\big(3^{-1}\|A\|^{-1}(5Ln^{-1}M)^{-\frac12},9^{-1}L^{-1}(15+7\ln M)^{-1}\big),\\
C_0&:=\big(14^2 \sqrt{L}M\|A\| \ln(2e^2n\sqrt{L}\|A\|^{-1})\big)^{-1}.
\end{align*}
Note that the constants $C_0$ and $\overline{C_0}$ are generally much smaller than $L^{-1}$, when compared with Assumption \ref{ass}(i).
\begin{theorem}\label{thm:main} 
Let Assumption \ref{ass} hold with $b=(1+2\nu)^{-1}$. 
Then there exists some $c^*$ independent of $k$, $n$ or $\delta$ such that, for any $k\geq 0$, 
\begin{align}\label{err:main1}
\E[\|e_k^\delta\|^2]^\frac12&\leq c^*k^{-\min(\nu,\frac12)}+c^*
\left\{\begin{array}{cc}
\delta^{\frac{2\nu}{1+2\nu}},     & a>0, \\
n^{-\frac12}\sqrt{k}\delta,     & a=0,
\end{array}\right.    & c_0<\overline{C_0}, \\
\label{err:main2}
\|e_k^\delta\|&\leq \sqrt{n}c^*k^{-\frac12+\max(\frac12-\nu,0)}+c^*
\left\{\begin{array}{cc}
\delta^{\frac{2\nu}{1+2\nu}},     & a>0, \\
n^{-\frac12}\sqrt{k}\delta,     & a=0,
\end{array}\right.    & c_0<C_0.
\end{align}
\end{theorem}

The next corollary follows directly from Theorem \ref{thm:main}. The shorthand notation $k(\delta)=\mathcal{O}(\delta^{-\frac{2}{1+2\nu}})$ means that there exist constants $0<C_1<C_2$ and some $\delta_0>0$, such that for any $\delta \in (0 ,\delta_0]$, there holds
$C_1\delta^{-\frac{2}{1+2\nu}} \leq k(\delta) \leq C_2\delta ^{- \frac{2}{1+2\nu}}.$
\begin{corollary}\label{cor:rate}
Under suitable step size schedules, the following statements hold.
\begin{enumerate}
\item[{\rm(i)}] When $a>0$ and $b=(1+2\nu)^{-1}$, i.e., rSVRG, there holds for any small $\epsilon>0$,
\begin{align*}
\E[\|e_k^\delta\|^2]^\frac12\leq& \mathcal{O}(\delta^{\frac{2\nu}{1+2\nu}}),\quad \forall k\geq k(\delta):= \delta^{-\frac{2\nu}{(1+2\nu)\min(\nu,\frac12)}},\\
\|e_k^\delta\|\leq& \mathcal{O}(\delta^{\frac{2\nu}{1+2\nu}}),\quad \forall k\geq k(\delta):= \delta^{-\frac{2\nu}{(1+2\nu)\min(\nu,\frac12-\epsilon)}}.
\end{align*}
\item[{\rm(ii)}] When $a=0$, i.e., SVRG, there holds
\begin{align*}
\E[\|e_{k(\delta)}^\delta\|^2]^\frac12\leq& \mathcal{O}(\delta^{\frac{2\nu}{1+2\nu}}),\quad \forall k(\delta)=\mathcal{O}(\delta^{-\frac{2}{1+2\nu}}),\; \nu\in(0,\tfrac12],\\ \|e_{k(\delta)}^\delta\|\leq & \mathcal{O}(\delta^{\frac{2\nu}{1+2\nu}}),\quad \forall k(\delta)=\mathcal{O}(\delta^{-\frac{2}{1+2\nu}}),\; \nu\in(0,\tfrac12).
\end{align*}
\end{enumerate}    
\end{corollary}

\begin{remark}\label{rem:rate}
Note that the convergence rate $\mathcal{O}(\delta^{\frac{2\nu}{1+2\nu}}) $ stated in Corollary \ref{cor:rate} is identical with that for the Landweber method \cite{EnglHankeNeubauer:1996}, which is known to be order-optimal. With a suitable choice of $A$, rSVRG can achieve optimal convergence rates without any early stopping rule, instead inheriting the regularizing property of the spectral cutoff. SVRG is also optimal with \textit{a priori} stopping rules for  $\nu\in(0,\frac12)$. 
These rates are identical with that of SVRG in \cite{JinZhouZou:2022ip,JinChen:2024}. Note that, when $\nu\in(0,\frac12]$, the condition for optimal convergence rates in expectation of standard SVRG is more relaxed than that in \cite{JinZhouZou:2022ip}, which requires also a special structure on $A_\dag$ and the step size $c_0\leq \mathcal{O}\big((Mn^{-1}\|A\|^2)^{-1}\big)$.
It is comparable with that in \cite{JinChen:2024} for small $M$ and more relaxed than that for relatively large $M$.
\end{remark}

\begin{remark}
Note that the error estimates in Corollary \ref{cor:rate} are derived under the assumption that the smoothness parameter $\nu$ is known \textit{a priori} 
for determining the stopping index $k(\delta)$ for SVRG and the truncation level $J$ in rSVRG. However, the parameter $\nu$ is generally unknown in practical applications. Thus the \textit{a priori} choice represents one significant restriction of the applicability of rSVRG. It is of much interest to develop \textit{a posteriori} or data-driven choices of the truncation level $J$.
\end{remark}

Assumption \ref{ass}(iii) is to simplify the analysis of rSVRG in Section \ref{sec:conv}. 
The result for rSVRG (i.e., $a>0$) remains valid when $A_\dag$ is approximated by some operator $A$ suitably; see the next corollary.
\begin{corollary}\label{cor:svd0}
Let Assumption \ref{ass}(i) and (ii) hold.  Suppose that either $A:X\rightarrow Y$ is invertible with $\|A^{-1}\|\leq a^{-1}\delta^{-b}$ or $A$ is compact with the nonzero singular value greater than $a\delta^b$ for some $a>0$ and $b>0$. If $\|A-A_\dag\|$ is sufficiently small, then Theorem \ref{thm:main} remains valid for (r)SVRG.   
\end{corollary}

In the absence of the source condition in Assumption \ref{ass}(ii), the regularizing property of (r)SVRG remains valid both in expectation and in the uniform sense.
\begin{corollary}\label{cor:regularizing}
Let Assumption \ref{ass}(i) and (iii) hold. Then rSVRG is regularizing $($in expectation when $c_0< \bar{C_0}$ and in the uniform sense when $c_0< C_0)$ when $b\in(0,1)$, and SVRG is regularizing $($in expectation when $c_0< \bar{C_0}$ and in the uniform sense when $c_0< C_0)$ when
the stopping index $k(\delta)$ satisfies the following conditions $$\lim_{\delta\to 0^+}k(\delta)=\infty \quad \mbox{and}\quad \lim_{\delta\to 0^+}\sqrt{k(\delta)}\delta=0.$$
\end{corollary}

\section{Convergence analysis}\label{sec:conv}

To prove Theorem \ref{thm:main}, we first give several shorthand notation. We denote the (r)SVRG iterates for the noisy data $y^\delta$ by $x_k^\delta$. For any $K=0,1,\cdots$ and $t=0,\cdots,M-1$, we define
\begin{align*}
&e_{KM+t}^\delta=x_{KM+t}^\delta-x_\dag,
\quad \Delta_{KM+t}^\delta=x_{KM+t}^\delta-x_{KM}^\delta,\\
&P_{KM+t}=I-c_0A_{i_{{KM+t}}}^*A_{i_{{KM+t}}}, \quad P=I-c_0B,\\
&N_{KM+t}=B-A_{i_{{KM+t}}}^*A_{i_{{KM+t}}}, \quad \mbox{and} \quad \zeta=n^{-1}A^* \xi,
\end{align*}
with $B:=\E[A_i^*A_i]=n^{-1}A^* A\;: X\rightarrow X$.
Then there hold
$$\Delta^\delta_{KM}=0, \quad \E[P_{KM+t}]=P \quad \mbox{and}\quad \E[N_{KM+t}]=0.$$
We also define the summations
\begin{align*}
\phi_{i}^{i'}&=\sum_{j=i}^{i'}P^{i'-j}N_j\Delta_j^\delta, \quad \tilde\phi^{i'}=\sum_{j=0}^{i'}P^{j}, \quad \forall i,\;i'\geq 0,\\
\overline{\Phi}_{i}^{i'}(0,1)&=2c_0\|B\|\E[\|A \Delta_{i'}^\delta\|^2]+\sum_{j=i}^{i'-1}(i'-j)^{-1}\E[\|A \Delta_j^\delta\|^2], \quad \forall i'\geq i\geq 0,\\
\Phi_{i}^{i'}(0,\tfrac12)&=\sqrt{2c_0\|B\|}\|A \Delta_{i'}^\delta\|+\sum_{j=i}^{i'-1}(i'-j)^{-\frac12}\|A \Delta_j^\delta\|, \quad \forall i'\geq i\geq 0,
\end{align*}
and for any $r,\;j'>0$ or $r=0$,
\begin{align*}
\overline{\Phi}_{i}^{i'}(j',r)&=\sum_{j=i}^{i'}(i'+j'-j)^{-r}\E[\|A \Delta_j^\delta\|^2], \\
\Phi_{i}^{i'}(j',r)&=\sum_{j=i}^{i'}(i'+j'-j)^{-r}\|A \Delta_j^\delta\|,
\end{align*}
and follow the conventions $\sum_{j=i}^{i'}R_j=0$ and $\sum_{j=-i}^{i'}R_j=\sum_{j=0}^{i'}R_j$ for any sequence $\{R_j\}_j$ and $0\leq i'< i$.
Under Assumption \ref{ass}(iii), $A_\delta:=A_\dag-A$ satisfies $A^* A_\delta=0$ and $ \|A_\delta\|< a\delta^{b}.$ Similarly, let $$B_\dag:=\E[A_{\dag,i}^* A_{\dag,i}]=n^{-1}A_\dag^* A_\dag\quad \mbox{and}\quad B_\delta:=B_\dag-B=n^{-1}A_{\delta}^* A_\delta.$$ Then $B^* B_\delta=0$ and $ \|B_\delta\|<n^{-1}a^2\delta^{2b}.$

\subsection{Error decomposition}
For any $K\geq 0$ and $0\leq t\leq M-1$, we  decompose the error $e_{KM+t+1}^\delta
\equiv x_{KM+t+1}^\delta-x_\dag$ and the weighted successive error $A\Delta_{KM+t}^\delta= 
A(x_{KM+t}^\delta-x_{KM}^\delta)$ between the $(KM+t)$th and $KM$th iterations  into the bias and variance components, which plays a crucial role in the subsequent analysis. (Note that the $KM$th iteration is an anchor point.)
\begin{lemma}\label{lem:bias-var}
Let Assumption \ref{ass}(i) hold. Then for any $K\geq 0$, $0\leq t\leq M-1$ and $k=KM+t+1$, there hold
\begin{align*}
&\E[e_{k}^\delta]=P^{k}e_0^\delta+c_0\tilde\phi^{k-1} \zeta, \quad
e_{k}^\delta-\E[e_{k}^\delta]=c_0\phi_{1}^{k-1},\\[1mm]
&\E[A\Delta_{KM+t}^\delta]=A(P^t-I)P^{KM}e_0^\delta+c_0 A P^{KM}\tilde\phi^{t-1} \zeta,\\[1mm]
&A\Delta_{KM+t}^\delta-A\E[\Delta_{KM+t}^\delta]=c_0A(P^t-I)\phi_{1}^{KM-1} +c_0A\phi_{KM+1}^{KM+t-1}.
\end{align*}
\end{lemma}
\begin{proof}
From the definitions of $e_j^\delta$, $\Delta_j^\delta$, $P$ and $N_j$, we derive
\begin{align*}
e_{KM+t+1}^\delta &=P_{KM+t}e_{KM+t}^\delta-c_0N_{KM+t}e_{KM}^\delta+c_0 \zeta
=P e_{KM+t}^\delta+c_0 N_{KM+t}\Delta_{KM+t}^\delta+c_0 \zeta\\[2mm]  
&=\ldots = P^{t+1}e_{KM}^\delta + c_0\sum_{j=1}^{t+1}P^{j-1}N_{KM+t+1-j}\Delta_{KM+t+1-j}^\delta+c_0 \tilde\phi^{t}\zeta\\
&=\ldots = P^{KM+t+1}e_{0}^\delta + c_0\phi_{0}^{KM+t}+c_0 \tilde\phi^{KM+t}\zeta.
\end{align*}
When $t=M-1$, this identity gives
\begin{align*}
e_{(K+1)M}^\delta &=P^{(K+1)M}e_{0}^\delta + c_0\phi_{0}^{(K+1)M-1}+c_0 \tilde\phi^{(K+1)M-1}\zeta.
\end{align*}
Then, with the convention $\sum_{j=i}^{i'}R_j=0$ for any sequence $\{R_j\}_j$ and $i'< i$, we have
\begin{align*}
A\Delta_{KM+t}^\delta=& A\Big(P^{KM+t}e_{0}^\delta + c_0\phi_{0}^{KM+t-1}+c_0 \tilde\phi^{KM+t-1}\zeta-P^{KM}e_{0}^\delta - c_0\phi_{0}^{KM-1}-c_0 \tilde\phi^{KM-1}\zeta\Big)\\
=&A(P^t-I)P^{KM}e_{0}^\delta +c_0 AP^{KM}\tilde\phi^{t-1}\zeta
+ c_0A(P^t-I)\phi_{0}^{KM-1}+ c_0A\phi_{KM}^{KM+t-1}.
\end{align*}
Finally, the identities $\E[N_j]=0$ and $\Delta_0^\delta=\Delta_{KM}^\delta=0$ imply the desired identities.
\end{proof}

Based on the triangle inequality, we bound the error $e_k^\delta$ by
\begin{align*}
\E[\|e_{k}^\delta\|^2]^\frac12\leq\|\E[e_{k}^\delta]\|+\E[\|e_{k}^\delta-\E[e_{k}^\delta]\|^2]^\frac12 \quad \mbox{and}\quad \|e_{k}^\delta\|\leq\|\E[e_{k}^\delta]\|+\|e_{k}^\delta-\E[e_{k}^\delta]\|.
\end{align*} 
The next lemma bounds the bias $\|\E[e_{k}^\delta]\|$ and variance $\E[\|e_{k}^\delta-\E[e_{k}^\delta]\|^2]$ (and $\|e_{k}^\delta-\E[e_{k}^\delta]\|$) in terms of the weighted successive error $\E[\|A\Delta_{j}^\delta\|^2]$ (and $\|A\Delta_{j}^\delta\|$), respectively.

\begin{lemma}\label{lem:decom_err}
Let Assumption \ref{ass}(i) hold. Then for any $k\geq 0$,
\begin{align*}
\|\E[e_{k}^\delta]\|\leq&\|P^{k}e_0^\delta\|+n^{-\frac12}c_0\|\tilde\phi^{k-1} B^{\frac12}\| \delta,\\
\E[\|e_{k}^\delta-\E[e_{k}^\delta]\|^2]\leq& \frac{c_0}{2}\overline{\Phi}_{1}^{k-1}(0,1)\quad \mbox{and}\quad
\|e_{k}^\delta-\E[e_{k}^\delta]\|\leq  \sqrt{\frac{nc_0}{2}}\Phi_{1}^{k-1}(0,\tfrac12).
\end{align*}
\end{lemma}
\begin{proof}
Lemma \ref{lem:bias-var} and the definitions $\zeta=n^{-1}A^* \xi$ and $B=n^{-1}A^*A$ yield
\begin{align*}
\|\E[e_{k}^\delta]\|\leq \|P^{k}e_0^\delta\|+n^{-1}c_0\|\tilde\phi^{k-1} A^* \xi\|
\leq \|P^{k}e_0^\delta\|+n^{-\frac12}c_0\|\tilde\phi^{k-1} B^{\frac12}\| \delta.
\end{align*}
Similarly, by Lemma \ref{lem:bias-var} and the identity $\E[\langle P^{k-1-i}N_i \Delta_i^\delta, P^{k-1-j}N_j \Delta_j^\delta\rangle|\mathcal{F}_j]=0$ for any $j>i$, we have
\begin{align*}
\E[\|e_{k}^\delta-\E[e_{k}^\delta]\|^2]
= &c_0^2\E[\|\phi_{1}^{k-1}\|^2]
= c_0^2\sum_{j=1}^{k-1} \E[\|P^{k-1-j} N_j\Delta_j^\delta\|^2].
\end{align*}
Then, by Lemmas \ref{lem:kernel} and \ref{lem:N}, we derive 
\begin{align*}
\E[\|e_{k}^\delta-\E[e_{k}^\delta]\|^2]
\leq c_0^2\sum_{j=1}^{k-1} \|P^{k-1-j} B^\frac12\|^2 \E[\|A\Delta_j^\delta\|^2]
\leq \frac{c_0}{2}\overline{\Phi}_{1}^{k-1}(0,1).
\end{align*}
Similarly, by the triangle inequality and Lemmas \ref{lem:kernel} and \ref{lem:N}, we obtain 
\begin{align*}
\|e_{k}^\delta-\E[e_{k}^\delta]\|
\leq&c_0\sum_{j=1}^{k-1}\| P^{k-1-j} N_j\Delta_j^\delta\|
\leq\sqrt{n}c_0\sum_{j=1}^{k-1}\| P^{k-1-j} B^\frac12\|\|A\Delta_j^\delta\|
\leq\sqrt{\frac{nc_0}{2}}\Phi_{1}^{k-1}(0,\tfrac12).
\end{align*}
This completes the proof of the lemma.
\end{proof}

Now we bound the weighted successive errors $\E[\|A\Delta_{j}^\delta\|^2]$ and $\|A\Delta_{j}^\delta\|$; see Appendix \ref{app:estimate} for the lengthy and technical proof.  
\begin{theorem}\label{thm:Delta}
Let Assumption \ref{ass}(i) hold. Then there exist some $c_1$ and $c_2$ independent of $k$, $n$, $\delta$ and $\nu$ such that, for any $k\geq 0$, 
\begin{align}\label{eqn:Delta_thmE1}
\E[\|A\Delta_{k}^\delta\|^2]&\leq (c_1+c_2\delta^2)(k+M)^{-2},\quad \mbox{if } c_0<\overline{C_0},\\
\label{eqn:Delta_thmE2}
\|A\Delta_{k}^\delta\|&\leq (c_1+c_2\delta)(k+M)^{-1},\quad \mbox{if } c_0<C_0.
\end{align}
\end{theorem}
\subsection{Convergence analysis}\label{sec:rates}

Now, using Theorem \ref{thm:Delta} and Lemma \ref{lem:decom_err}, we can prove Theorem \ref{thm:main}.

\begin{proof}
For any $k\geq 1$, the triangle inequality and Lemma \ref{lem:decom_err} give 
\begin{align}\label{eqn:err0}
\E[\|e_{k}^\delta\|^2]^\frac12\leq &\|P^{k}e_0^\delta\|+n^{-\frac12}c_0\|\tilde\phi^{k-1} B^{\frac12}\| \delta+\sqrt{\frac{c_0}{2}\overline{\Phi}_{1}^{k-1}(0,1)}.
\end{align}
When $c_0<\overline{C_0}$, by Theorem \ref{thm:Delta}, the estimate \eqref{eqn:sum01_e} in Lemma \ref{lem:sum} implies
\begin{align*}
\overline{\Phi}_{1}^{k-1}(0,1)
\leq (3+2c_0\|B\|)(c_1+c_2\delta^2)k^{-1}.
\end{align*}
Next, we bound the first two terms in \eqref{eqn:err0}. 
By the definitions $B_\dag=B+B_\delta$, $P=I-c_0B$ and $B=n^{-1}A^*A$ with $A$ satisfying Assumption \ref{ass}(iii), Assumption \ref{ass}(ii) implies 
\begin{equation*} 
e_0^\delta=B_\dag^\nu w=(B+B_\delta)^\nu w=B^\nu w+B_\delta^\nu w\quad \mbox{and}\quad 
PB_\delta^\nu w=B_\delta^\nu w.
\end{equation*}
These identities follow directly from the definition of the fractional power $B_\dag^\nu w$ via spectral decomposition:
\begin{align}
   & B_\dag^\nu w  = \sum_{j=1}^\infty (n^{-1}\sigma_j^2)^\nu \langle \varphi_j,w\rangle \varphi_j \nonumber\\
    =& \sum_{j=1}^J (n^{-1}\sigma_j^2)^\nu \langle \varphi_j,w\rangle\varphi_j + \sum_{j=J+1}^\infty (n^{-1}\sigma_j^2)^\nu \langle \varphi_j,w\rangle\varphi_j = B^\nu w + B_\delta^\nu w. \label{eqn:B^nu}
\end{align}
Together with Lemma \ref{lem:kernel} and the estimate $\|B_\delta\|<n^{-1}a^2\delta^{2b}$, we obtain 
\begin{align}
\|P^{k}e_0^\delta\|=&\|P^{k}B^\nu w+P^{k}B_\delta ^\nu w\|=\|P^{k}B^\nu w+ B_\delta ^\nu w\| \nonumber \\
\leq& (\|P^{k}B^\nu\|+\|B_\delta\|^{\nu}) \|w\|
\leq  (\nu^{\nu} c_0^{-\nu}k^{-\nu}+n^{-\nu}a^{2\nu}\delta^{2b\nu}) \|w\|.\label{eqn:approx}
\end{align}
Next we bound the term ${\rm I}:=n^{-\frac12}c_0\|\tilde\phi^{k-1} B^{\frac12}\|\delta$.
If $a=0$, Lemma \ref{lem:kernel}, the triangle inequality and inequality $c_0\|B\|\leq \max(\overline{C_0},C_0)L\leq \frac12$ imply 
\begin{align}\label{eqn:noise0}
{\rm I}
\leq& n^{-\frac12}c_0\sum_{j=0}^{k-1}\|P^j B^{\frac12}\|\delta
\leq \sqrt{\frac{c_0}{2n}}\Big(\sqrt{2c_0\|B\|}+\sum_{j=1}^{k-1}j^{-\frac12}\Big)\delta\nonumber\\
\leq& \sqrt{\frac{c_0}{2n}}\Big(1+1+\int_1^{k}x^{-\frac12}\; {\rm d}x\Big)\delta \nonumber\\
\leq& \sqrt{\frac{c_0}{2n}}(2+2\sqrt{k}-2)\delta\leq \sqrt{\frac{2c_0}{n}}\sqrt{k}\delta;
\end{align}
if $a>0$, for any $\lambda$ in the spectrum $\mathrm{Sp}(B)$ of $B$, either $\lambda\geq n^{-1}a^2 \delta^{2b}$ or $\lambda=0$ holds, and thus
\begin{align}
{\rm I}
&\leq n^{-\frac12}c_0\delta\sup_{\lambda\in \mathrm{Sp}(B)}\sum_{j=0}^{k-1}(1-c_0 \lambda)^j \lambda^{\frac12}\nonumber
\leq n^{-\frac12}c_0\delta\sup_{\lambda\geq n^{-1}a^2 \delta^{2b}}\big(1-(1-c_0 \lambda)^k\big)c_0^{-1}\lambda^{-\frac12}\\
&\leq n^{-\frac12}\delta\sup_{\lambda\geq n^{-1}a^2 \delta^{2b}}\lambda^{-\frac12}
\leq a^{-1}\delta^{1-b}.\label{eqn:noise1}
\end{align}
Since $b=(1+2\nu)^{-1}$ and $\delta<1$, we derive from \eqref{eqn:err0} and the above estimates that, when $a=0$, 
\begin{align*}
\E[\|e_{k}^\delta\|^2]^\frac12 \leq &\nu^{\nu} c_0^{-\nu}k^{-\nu} \|w\|+\sqrt{2c_0}n^{-\frac12}\sqrt{k}\delta+\sqrt{2c_0}(\sqrt{c_1}+\sqrt{c_2}\delta)k^{-\frac12}\\
\leq &\big(\nu^{\nu} c_0^{-\nu} \|w\|+\sqrt{2c_0}(\sqrt{c_1}+\sqrt{c_2})\big)k^{-\min(\nu,\frac12)}+\sqrt{2c_0}n^{-\frac12}\sqrt{k}\delta;
\end{align*}
and when $a>0$, 
\begin{align*}
\E[\|e_{k}^\delta\|^2]^\frac12\leq &(\nu^{\nu} c_0^{-\nu}k^{-\nu}+n^{-\nu}a^{2\nu}\delta^{2b\nu}) \|w\|+a^{-1}\delta^{1-b}+\sqrt{2c_0}(\sqrt{c_1}+\sqrt{c_2}\delta)k^{-\frac12}\\
\leq &\big(\nu^{\nu} c_0^{-\nu}\|w\|+\sqrt{2c_0}(\sqrt{c_1}+\sqrt{c_2})\big)k^{-\min(\nu,\frac12)} 
+(n^{-\nu}a^{2\nu}\|w\|+a^{-1})\delta^{\frac{2\nu}{1+2\nu}}.
\end{align*}
This proves the estimate \eqref{err:main1}.
Similarly, for $\|e_{k}^\delta\|$, when $c_0<C_0$, Lemma \ref{lem:decom_err} yields 
\begin{align}\label{eqn:err0_as}
\|e_{k}^\delta\|\leq &\|P^{k}e_0^\delta\|+n^{-\frac12}c_0\|\tilde\phi^{k-1} B^{\frac12}\| \delta
+\sqrt{\frac{nc_0}{2}}\Phi_{1}^{k-1}(0,\tfrac12).
\end{align}
Theorem \ref{thm:Delta} and the inequality \eqref{eqn:sum00.5_e} in Lemma \ref{lem:sum} imply
\begin{align*}
\Phi_{1}^{k-1}(0,\tfrac12)
\leq 3\sqrt{2}(c_1+c_2\delta)k^{-\frac12}\ln k.
\end{align*}
Then, by the conditions $b=(1+2\nu)^{-1}$ and $\delta<1$, we derive from \eqref{eqn:err0_as} and the estimates \eqref{eqn:approx}--\eqref{eqn:noise1} that, when $a=0$, 
\begin{align*}
\|e_{k}^\delta\|\leq &\nu^{\nu} c_0^{-\nu} \|w\|k^{-\nu}
+\sqrt{2c_0}n^{-\frac12}\sqrt{k}\delta
+3\sqrt{nc_0}(c_1+c_2\delta)k^{-\frac12}\ln k\\
\leq &\big(\nu^{\nu} c_0^{-\nu} \|w\|+3\sqrt{nc_0}(c_1+c_2)\big)k^{-\frac12+\max(\frac12-\nu,0)}
+\sqrt{2c_0}n^{-\frac12}\sqrt{k}\delta;
\end{align*}
and when $a>0$, 
\begin{align*}
\|e_{k}^\delta\|\leq &(\nu^{\nu} c_0^{-\nu}k^{-\nu}+n^{-\nu}a^{2\nu}\delta^{2b\nu}) \|w\|
+a^{-1}\delta^{1-b}
+3\sqrt{nc_0}(c_1+c_2\delta)k^{-\frac12}\ln k\\
\leq &\big(\nu^{\nu} c_0^{-\nu}\|w\|+3\sqrt{nc_0}(c_1+c_2)\big)k^{-\frac12+\max(\frac12-\nu,0)}+\big(n^{-\nu}a^{2\nu}\|w\|+a^{-1}\big)\delta^{\frac{2\nu}{1+2\nu}}.
\end{align*}
This proves the estimate \eqref{err:main2}, and completes the proof of the theorem.
\end{proof}

\begin{remark}\label{rem:a}
The parameter $b$ in Assumption \ref{ass}{\rm(}iii{\rm)} is set to $b=(1+2\nu)^{-1}$. Now we discuss the suitable choice of the parameter $a$. 
From the bound on $\|e_{k}^\delta\|$ {\rm(}or $\E[\|e_{k}^\delta\|^2]^\frac12${\rm)} in the proof of Theorem \ref{thm:main}, we have
\begin{align*}
\lim_{k\to \infty} \|e_{k}^\delta\|\leq c_{n,\nu}(a)\delta^{\frac{2\nu}{1+2\nu}}, \quad \mbox{with }c_{n,\nu}(a)=n^{-\nu}a^{2\nu}\|w\|+a^{-1}.
\end{align*}
Note that $c_{n,\nu}$ attains its minimum at $a_*:=\big(n^\nu/(2\nu\|w\|)\big)^{\frac{1}{1+2\nu}}$, with 
$$c_{n,\nu}(a_*)=\Big((2\nu)^{-\frac{2\nu}{1+2\nu}}+(2\nu)^{\frac{1}{1+2\nu}}\Big)n^{-\frac{\nu}{1+2\nu}}\|w\|^{\frac{1}{1+2\nu}}.$$
To avoid the blow-up of $a$ as $\nu\to 0^+$, let $a=(n^\nu/\|w\|)^{\frac{1}{1+2\nu}}$. Then $c_{n,\nu}(a)= 2n^{-\frac{\nu}{1+2\nu}} \|w\|^{\frac{1}{1+2\nu}}$.
\end{remark}

Next we prove Corollary \ref{cor:svd0}, which relaxes Assumption \ref{ass}(iii).
\begin{proof}
When $a>0$, let $\|A_\dag-A\|\leq \epsilon_A\leq \|A_\dag\|$ and $B=n^{-1} A^* A$. 
Under Assumption \ref{ass}(ii), we can bound the term $\|P^k e_0^\delta\|$ in \eqref{eqn:err0} and \eqref{eqn:err0_as} by
\begin{align*}
\|P^k e_0^\delta\|=&\|P^k B_\dag^\nu w\|
\leq \|P^k B^\nu w\| +\|P^k (B_\dag^\nu-B^\nu)w\|\\
\leq& \|P^k B^\nu w\| +\|B_\dag^\nu-B^\nu\|\|w\|.
\end{align*}
When $\nu\in(0,1]$, by \cite[Theorem 2.3]{KittanehKosaki:1987}, the term ${\rm I}:=\|B_\dag^\nu-B^\nu\|$ can be bounded by 
\begin{align*}
{\rm I} \leq& \|B_\dag-B\|^\nu = n^{-\nu}\|A_\dag^* A_\dag-A^* A\|^\nu\\
\leq &n^{-\nu}\big(\|A_\dag^*\| \|A_\dag-A\|+\|A_\dag^*-A^*\| \|A\|\big)^\nu\\
\leq& n^{-\nu}\big(2\|A_\dag\|+\epsilon_A\big)^\nu \epsilon_A^\nu
\leq \big(3n^{-1} \|A_\dag\|\big)^\nu \epsilon_A^\nu. 
\end{align*}
When $\nu=1$, $\|B_\dag-B\|\leq 3n^{-1}\|A_\dag\|\epsilon_A$.
When $\nu>1$, the function $h(z):=z^\nu$ is Lipchitz continuous on any closed interval in $[0,\infty)$, and thus
\begin{align*}
{\rm I} \leq& \nu \max\big(\|B_\dag\|,\|B\|\big)^{\nu-1}\|B_\dag-B\|\\
\leq& n^{-(\nu-1)}\nu \big(\|A_\dag\|+\epsilon_A\big)^{2(\nu-1)}\|B_\dag-B\|\\
\leq &2^{2\nu}\nu n^{-\nu} \|A_\dag\|^{2\nu-1}\epsilon_A.
\end{align*}
Then, let $\epsilon_A\leq \delta^{\frac{2\nu}{(1+2\nu)\min(1,\nu)}}$, we have 
$${\rm I}\leq c(\nu)n^{-\nu}\epsilon_A \leq c(\nu)n^{-\nu} \delta^{\frac{2\nu}{1+2\nu}},$$ 
with the constant $c(\nu)$ independent of $\delta$ and $n$. 
The assumption on $A$ implies $\|A^{-1}\|\leq a^{-1}\delta^{-b}$ or the nonzero singular values $\sigma$ of $A$ such that $\sigma\geq a\delta^b>0$, which implies \eqref{eqn:noise1}. Thus, Theorem \ref{thm:main} still holds. 
\end{proof}

The next remark complements Corollary \ref{cor:svd0} when $A_\dag$ is compact
and has an approximate truncated SVD $A$. \begin{remark}
\noindent If $A_\dag$ is compact, with its SVD $A_\dag(\cdot)=\sum_{j=1}^\infty \sigma_{j}\langle\varphi_{j},\cdot\rangle\psi_{j}$, where the singular values $\{\sigma_j\}_{j=1}^\infty$ such that $\sigma_j\geq\sigma_{j'}\geq 2a\delta^b>0$ for any $j\leq j'\leq J$ and $\sigma_j<2a\delta^b$ for any $j>J$.
For any small $\epsilon_A\in(0, a\delta^b)\subset(0,\|A^\dag\|)$, we may approximate the operator $A_\dag$ by $A(\cdot)=\sum_{j=1}^{J} \tilde{\sigma}_{j}\langle\tilde{\varphi}_{j},\cdot\rangle\tilde{\psi}_{j}$ with $\{\tilde{\varphi}_{j}\}_{j=1}^{J}$ and $\{\tilde{\psi}_{j}\}_{j=1}^{J}$ being orthonormal in $X$ and $Y$, respectively, which satisfies $$\|\varphi_{j}-\tilde{\varphi}_{j}\|<\epsilon_{A}\quad\mbox{and}\quad |\tilde\sigma_{j}-\sigma_{j}|<\epsilon_{A}.$$ Then we take $A_{\rm T}(\cdot)=\sum_{j=1}^{J} \sigma_{j}\langle\varphi_{j},\cdot\rangle\psi_{j}$. Let 
\begin{align*} 
B(\cdot)&=n^{-1} A^* A(\cdot)=n^{-1} \sum_{j=1}^{J} \tilde{\sigma}_{j}^2\langle\tilde{\varphi}_{j},\cdot\rangle\tilde{\varphi}_{j},\\ B_{\rm T}(\cdot)&=n^{-1} A_{\rm T}^* A_{\rm T}(\cdot)=n^{-1} \sum_{j=1}^{J} \sigma_{j}^2\langle \varphi_{j},\cdot\rangle \varphi_{j},
\end{align*}
and $B_{\rm T,\delta}=B_\dag-B_{\rm T}$. Then there hold 
$B_{\rm T}^* B_{\rm T,\delta}=0$, $\|B_{\rm T,\delta}\|<4n^{-1}a^2\delta^{2b}$, and 
$B_\dag^\nu = B_{\rm T}^\nu+B_{\rm T,\delta}^\nu$ due to the identity \eqref{eqn:B^nu}.
Hence,
\begin{align*}
{\rm I} &= \|B_\dag^\nu-B^\nu\|
= \|B_{\rm T}^\nu+B_{\rm T,\delta}^\nu-B^\nu\|\leq \|B_{\rm T,\delta}\|^\nu+{\rm I}_\epsilon, 
\end{align*}
with 
$$ {\rm I}_{\epsilon}:=n^{-\nu}\sup_{\|z\|=1}\big\|  \sum_{j=1}^{J} \big(\sigma_{j}^{2\nu}\langle \varphi_{j},z\rangle \varphi_{j}- \tilde{\sigma}_{j}^{2\nu}\langle\tilde{\varphi}_{j},z\rangle\tilde{\varphi}_{j}\big)\big\|.$$ 
By the triangle inequality, we obtain
\begin{align*}
{\rm I}_{\epsilon}
\leq&  n^{-\nu}\sup_{\|z\|=1}\bigg\|  \sum_{j=1}^{J} \sigma_{j}^{2\nu}\langle \varphi_{j},z\rangle (\varphi_{j}- \tilde{\varphi}_{j})\bigg\|\\
&+n^{-\nu}\sup_{\|z\|=1}\bigg\|  \sum_{j=1}^{J} (\sigma_{j}^{2\nu}- \tilde{\sigma}_{j}^{2\nu})\langle \varphi_{j},z\rangle\tilde{\varphi}_{j}\bigg\|\\
&+n^{-\nu}\sup_{\|z\|=1}\bigg\|  \sum_{j=1}^{J} \tilde{\sigma}_{j}^{2\nu}\langle\varphi_{j}-\tilde{\varphi}_{j},z\rangle\tilde{\varphi}_{j}\bigg\| \\
\leq& n^{-\nu}(\|A_\dag\|^{2\nu}+\|A\|^{2\nu})  \epsilon_A+n^{-\nu}\sup_{j\leq J}{|\sigma_j^{2\nu}-\tilde\sigma_j^{2\nu}|} \\
\leq& n^{-\nu}(1+2^{2\nu})\|A_\dag\|^{2\nu}  \epsilon_A+2\nu n^{-\nu} \sup_{j\leq J}\max(\sigma_j^{2\nu-1},\tilde\sigma_j^{2\nu-1})\epsilon_A\\
\leq& n^{-\nu}(1+2^{2\nu})\|A^\dag\|^{2\nu}  \epsilon_A+2\nu n^{-\nu} \max\big(( a\delta^b)^{2\nu-1},(2 \|A_\dag\|)^{2\nu-1}\big)\epsilon_A.
\end{align*}
Let $b=(1+2\nu)^{-1}$. If $\epsilon_A\leq \delta^{\frac{\max(1,2\nu)}{1+2\nu}}$, then ${\rm I}\leq \tilde{c}(\nu) n^{-\nu} \delta^{\frac{2\nu}{1+2\nu}}$, with  $\tilde{c}(\nu)$ independent of $\delta$ and $n$. 
The condition $\tilde{\sigma}_j\geq \sigma_j-\epsilon_A\geq a\delta^b>0$ for any $j\leq J$ implies \eqref{eqn:noise1}, and Theorem \ref{thm:main} still holds.
\end{remark}

Last, we give the proof of Corollary \ref{cor:regularizing}.

\begin{proof}
Note that the initial error $x_\dag-x_0\in \overline{\mathrm{Range}(A_\dag^*)}$. The polar decomposition $A_\dag=Q(A_\dag^* A_\dag)^\frac12$ with a partial isometry $Q$ (i.e., $Q^* Q$ and $Q Q^*$ are projections) implies $x_\dag-x_0\in \overline{\mathrm{Range}\big((A_\dag^* A_\dag)^\frac12\big)}$.
Thus, for any $\epsilon_0>0$, there exists some $\tilde{x}_0$, satisfying Assumption \ref{ass}(ii) with $\nu=\frac12$ and $w\equiv w_{\epsilon_0}\in Y$ (which depends on $\epsilon_0$), such that $\|x_0-\tilde{x}_0\|<\epsilon_0$.  
Let $\tilde{x}^\delta_k$ be the $k$th (r)SVRG iterate starting with $\tilde{x}_0$ and $\tilde{e}^\delta_k=\tilde{x}^\delta_k-x_\dag$. 
Then, for any $b\in(0,1)$, by Lemma \ref{lem:kernel} and the inequality \eqref{eqn:approx}, we can bound  $\|P^k e_0^\delta\|$ in \eqref{eqn:err0} and \eqref{eqn:err0_as} by
\begin{align*}
\|P^k e_0^\delta\|&\leq \|P^k \tilde{e}_0^\delta\|+\|P^k(\tilde{e}_0^\delta-e_0^\delta)\|\leq  \|P^kB_\dag^{\frac12}w\|+\epsilon_0\\
&\leq \big((2c_0)^{-\frac12}k^{-\frac12}+n^{-\frac12}a \delta^b\big) \|w\| +\epsilon_0.
\end{align*}
Consequently, for some constant $c^*$ depending on $\epsilon_0$, there hold
\begin{align*}
\begin{aligned}
\E[\|e_k^\delta\|^2]^\frac12&\leq \epsilon_0+c^*k^{-\frac12}+c^*
\left\{\begin{array}{cc}
\delta^{\min(b,1-b)},     & a>0, \\
n^{-\frac12}\sqrt{k}\delta,     & a=0,
\end{array}\right.    \quad c_0<\overline{C_0},\\
\|e_k^\delta\|&\leq \epsilon_0+\sqrt{n}c^*k^{-\frac12}\ln k+c^*
\left\{\begin{array}{cc}
\delta^{\min(b,1-b)}, & a>0, \\
n^{-\frac12}\sqrt{k}\delta,     & a=0,
\end{array}\right.    \quad c_0<C_0.
\end{aligned}
\end{align*}
Now for any fixed $\varepsilon>0$, there exists some $\tilde{x}_0$ such that $\epsilon_0<\frac13 \varepsilon$. 
Given the constant $c^*$ (depending on $\epsilon_0$), there exists some $k_1\in \mathbb{N}$ such that for any $k\geq k_1$, we have $\sqrt{n}c^*k^{-\frac12}\ln k<\frac13 \varepsilon$. Now we discuss the two cases (i) $a>0$ and $b\in (0,1)$, and (ii) $a=0$ separately.

\noindent \underline{Case (i) $a>0$ and $b\in(0,1)$.} Then there exists some $\delta_1>0$ such that for any $\delta\leq \delta_1$, there holds $c^* \delta^{\min(b,1-b)}<\frac13 \varepsilon$, which implies that for any $k\geq k_1$, there hold
\begin{align}
\label{eqn:err-limit1}\lim_{\delta\to0^+}\E[\|e_{k}^\delta\|^2]^\frac12&<\varepsilon,\quad \mbox{if }c_0<\overline{C_0},\\
\label{eqn:err-limit2}\lim_{\delta\to0^+}\|e_{k}^\delta\|&<\varepsilon,\quad \mbox{if }c_0<C_0.
\end{align}

\noindent\underline{Case (ii) $a=0$. }
Similar to case (i), when $a=0$, there exists some $\delta_2>0$ such that for any $\delta\leq \delta_2$, there holds $c^* n^{-\frac12}\sqrt{k_1}\delta<\frac13 \varepsilon$. 
Further, for any fixed $\delta\leq \delta_2$, there exists some $k_2(\delta)\geq k_1$ such that for any $k_1\leq k(\delta)\leq k_2(\delta)$, there hold 
$$\sqrt{n}c^* k(\delta)^{-\frac12}\ln k(\delta)<\tfrac13 \varepsilon\quad \mbox{and}\quad c^* n^{-\frac12}\sqrt{k(\delta)}\delta<\tfrac13 \varepsilon.$$
This discussion implies that for the stopping index $k(\delta)$ satisfying the following conditions $$\lim_{\delta\to 0^+}k(\delta)=\infty \quad \mbox{and}\quad \lim_{\delta\to 0^+}\sqrt{k(\delta)}\delta=0,$$ there again hold the relations \eqref{eqn:err-limit1} and \eqref{eqn:err-limit2}.\\ 
Combining the preceding discussions completes the proof of the corollary. 
\end{proof}

\section{Numerical experiments and discussions}\label{sec:num}
In this section, we provide numerical experiments for several linear inverse problems to complement the theoretical findings in Section \ref{sec:main}. The experimental setting is identical to that in \cite{JinZhouZou:2022ip}. We employ three examples,  i.e., \texttt{s-phillips} (mildly ill-posed), \texttt{s-gravity} (severely ill-posed) and \texttt{s-shaw}
(severely ill-posed), which are generated from the code \texttt{phillips}, \texttt{gravity} and \texttt{shaw}, taken from the \texttt{MATLAB}
package Regutools \cite{P.C.Hansen2007} (publicly available at \url{http://people.compute.dtu.dk/pcha/Regutools/}). All the examples are discretized into a finite-dimensional linear system with the forward operator $A_\dag: \mathbb{R}^m \rightarrow \mathbb{R}^n$ of size $n = m = 1000$, with $A_\dag x= (A_{\dag,1} x,\cdots,A_{\dag,n} x)$ for all $x\in \mathbb{R}^m$ and $A_{\dag,i}: \mathbb{R}^m \rightarrow \mathbb{R}$. 
To precisely control the regularity index $\nu$ in the source condition in Assumption
\ref{ass}(ii), we generate the exact solution $x_\dag$ by 
\begin{equation}\label{eqn:num_gene_x} 
x_\dag = \|(A_\dag^*A_\dag)^\nu x_e\|_{\ell^\infty}^{-1}(A_\dag^*A_\dag)^\nu x_e,
\end{equation}
with $x_e$ being the exact solution provided by the package and $\|\cdot\|_{\ell^\infty}$ the maximum norm of a vector. Note that the smoothness index $\nu$ in the source condition is slightly larger than the one in \eqref{eqn:num_gene_x} due to the existing regularity $\nu_e$ of $x_e$. The exact data $y_\dag$ is given by $y_\dag=A_\dag x_\dag$ and the noisy data $y^\delta$ is generated by
$$y^\delta_i:=y_{\dag,i}+\epsilon\|y_\dag\|_{\ell^\infty}\xi_i, \quad i=1,\cdots,n,$$
where the noise components $\xi_i$s follow the standard normal distribution, and $\epsilon > 0$ is the relative noise level. 

All the iterative methods are initialized to zero, with a constant step size $c_0=\|A_\dag\|^{-2}$ for the Landweber method (LM) and $c_0=\mathcal{O}(c)$ for (r)SVRG and SGD, where $c=\min_i(\|A_i\|^{-2})=L^{-1}$. 
The constant step size $c_0$ is taken for rSVRG so as to achieve optimal convergence while maintaining computational efficiency across all noise levels.  The methods are run for a maximum 1e5 epochs, where one epoch refers to one Landweber iteration, $nM/(n+M)$ (r)SVRG iterations, or $n$ SGD iterations, so that their overall computational complexity is comparable. 
The frequency $M$ of computing the full gradient is set to $M=2n$ as suggested in \cite{JohnsonZhang:2013}.
The operator $A$ for rSVRG is generated by the truncated SVD of $A_\dag$ with $b=1/\big(1+2(\nu+\nu_e)\big)$ and $a=(\|A_\dag\|/\|y_\dag\|)(n^{\nu+\nu_e}/c_1)^{\frac{1}{1+2(\nu+\nu_e)}}$, cf. Theorem \ref{thm:main} and Remark \ref{rem:a}, in which $\nu_e$ is taken to be 0.25 for \texttt{phillips}, 0.1 for \texttt{gravity}, and 0 for \texttt{shaw}. Note the constant $c_1$ is fixed for each problem with different regularity indices $\nu$ and noise levels $\epsilon$. 
One can also use the randomized SVD to generate $A$. Note that the \textit{a priori} choice of the parameters $a$ and $b$ is of theoretical interest only due to their dependence on the smoothness parameter $\nu+\nu_e$, which however is often unknown in practice.

For LM, the stopping index $k_*=k(\delta)$ (measured in terms of epoch count) is chosen by the discrepancy principle with $\tau =1.01$:
\begin{align*}
k(\delta)=\min\{k\in \mathbb{N}:\;\|A_\dag x_{k}^\delta-y^\delta\|\leq \tau \delta\},
\end{align*}
which can achieve order optimality. For rSVRG, $k_*$ is selected to be greater than the last index at which the iteration error exceeds that of LM upon its termination or the first index for which the iteration trajectory has plateaued. 
For SVRG and SGD, $k_*$ is taken such that the error is the smallest along the iteration trajectory.
The accuracy of the reconstructions is measured by the relative error $e_*={\E[\| x_{k_*}^\delta-x_\dag\|^2]^\frac12}/{\|x_\dag\|}$ for both (r)SVRG and SGD, 
and $e_*=\|x_{k_*}^\delta-x_\dag\|/\|x_\dag\|$ for LM. 
The statistical quantities generated by (r)SVRG are computed based on ten independent runs. 

\begin{figure}[hbt!]
\centering
  \setlength{\tabcolsep}{4pt}
\begin{tabular}{ccc}
\includegraphics[width=0.31\textwidth,trim={1.5cm 0 0.5cm 0.5cm}]{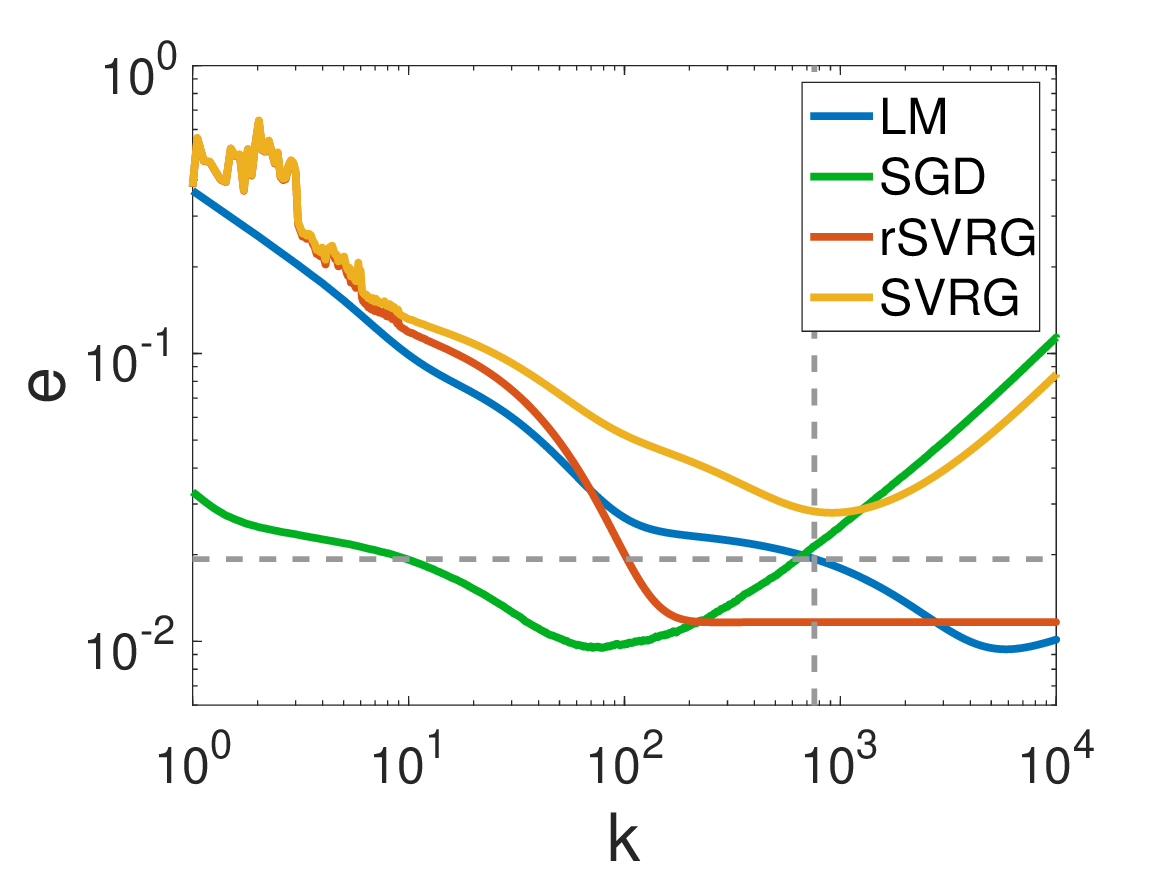}&
\includegraphics[width=0.31\textwidth,trim={1.5cm 0 0.5cm 0.5cm}]{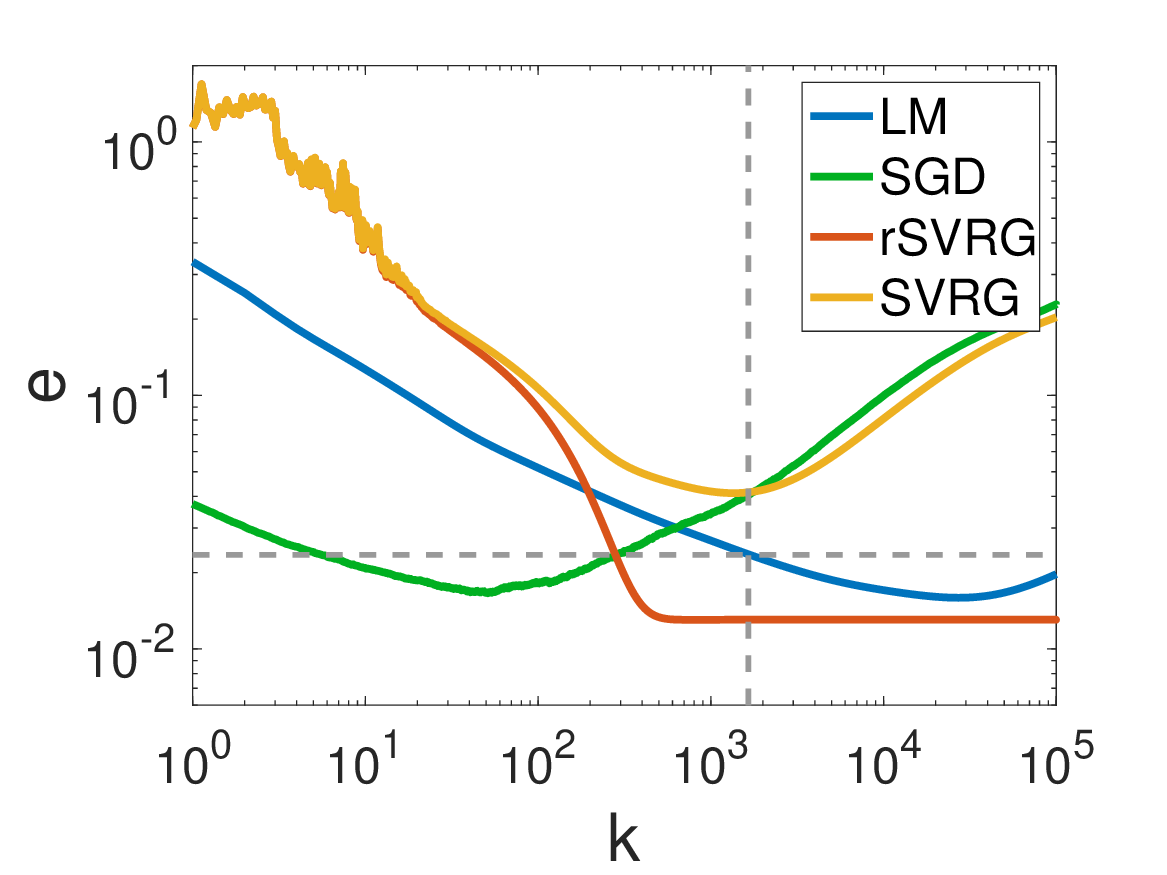}&
\includegraphics[width=0.31\textwidth,trim={1.5cm 0 0.5cm 0.5cm}]{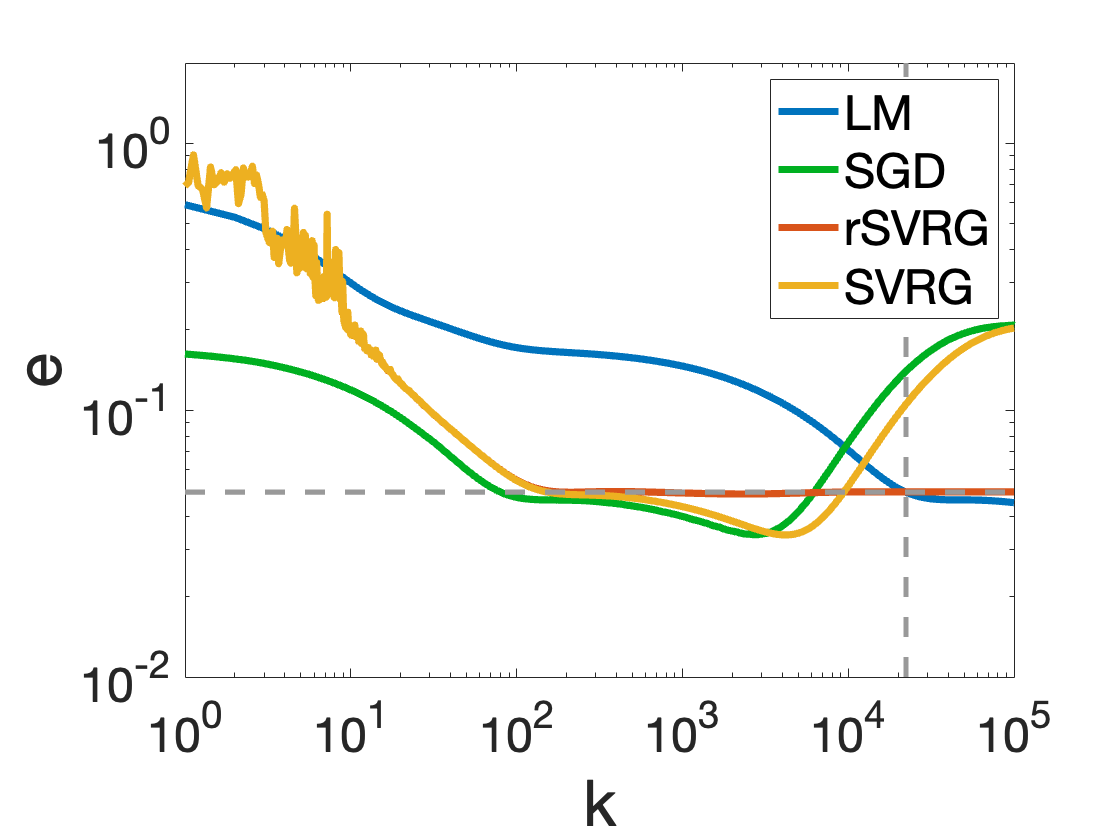}\\
\includegraphics[width=0.31\textwidth,trim={1.5cm 0 0.5cm 0.5cm}]{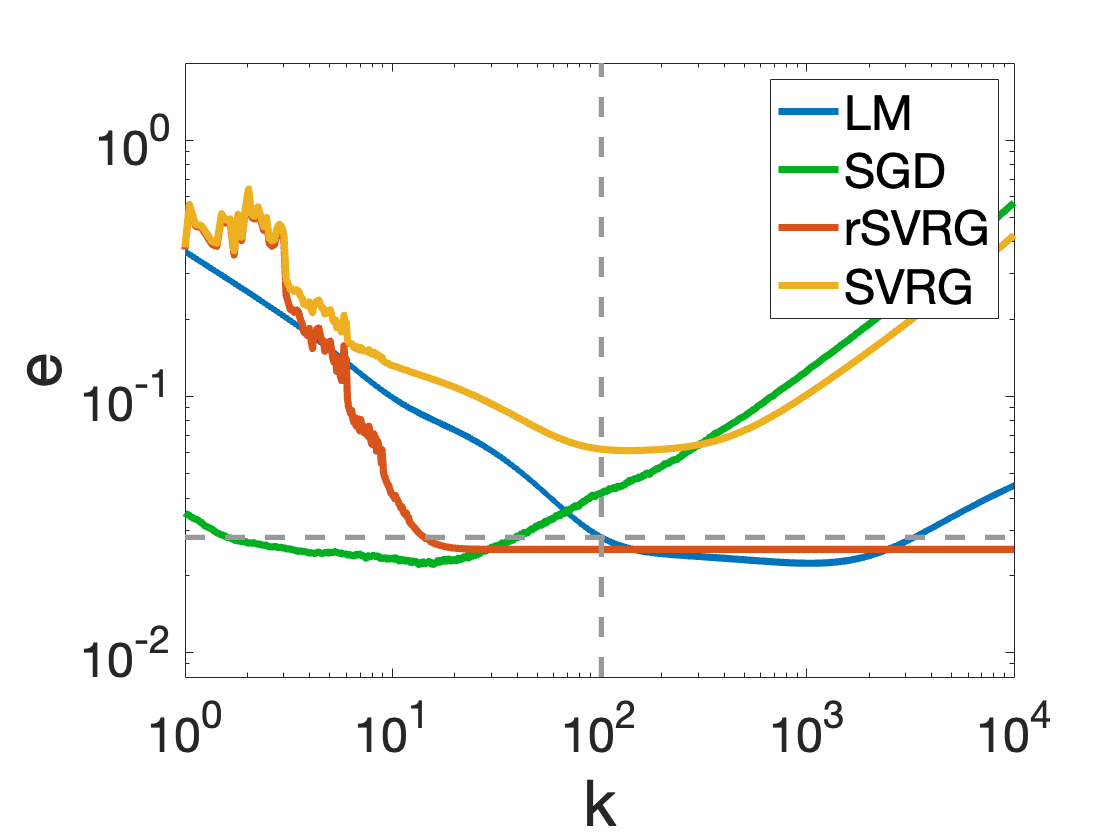}&
\includegraphics[width=0.31\textwidth,trim={1.5cm 0 0.5cm 0.5cm}]{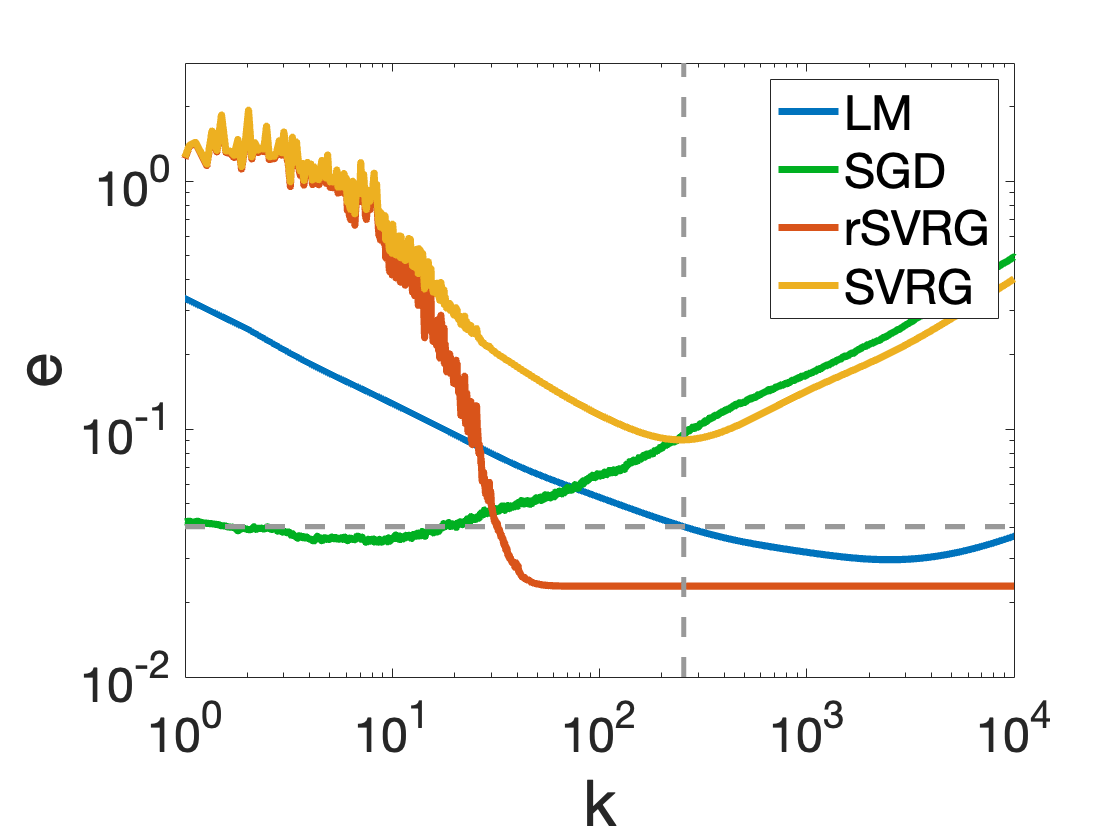}&
\includegraphics[width=0.31\textwidth,trim={1.5cm 0 0.5cm 0.5cm}]{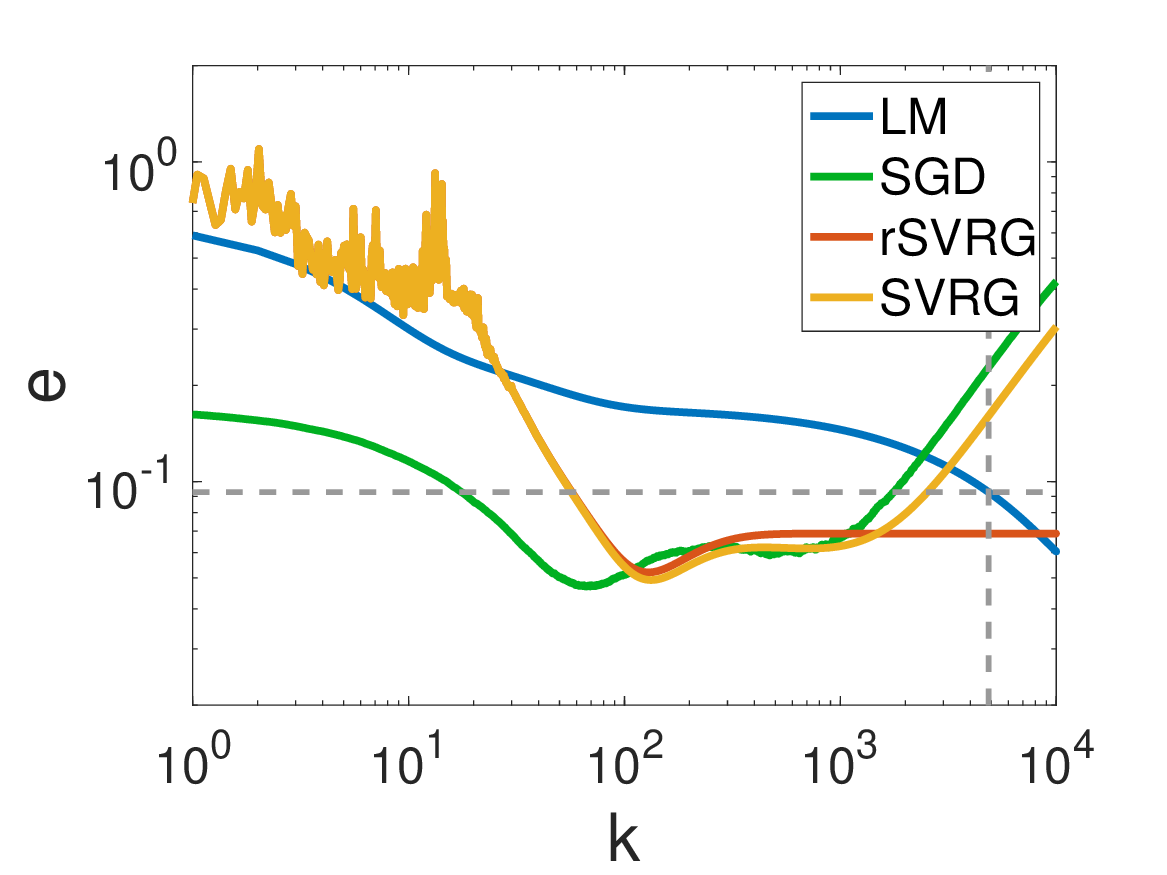}\\
\includegraphics[width=0.31\textwidth,trim={1.5cm 0 0.5cm 0.5cm}]{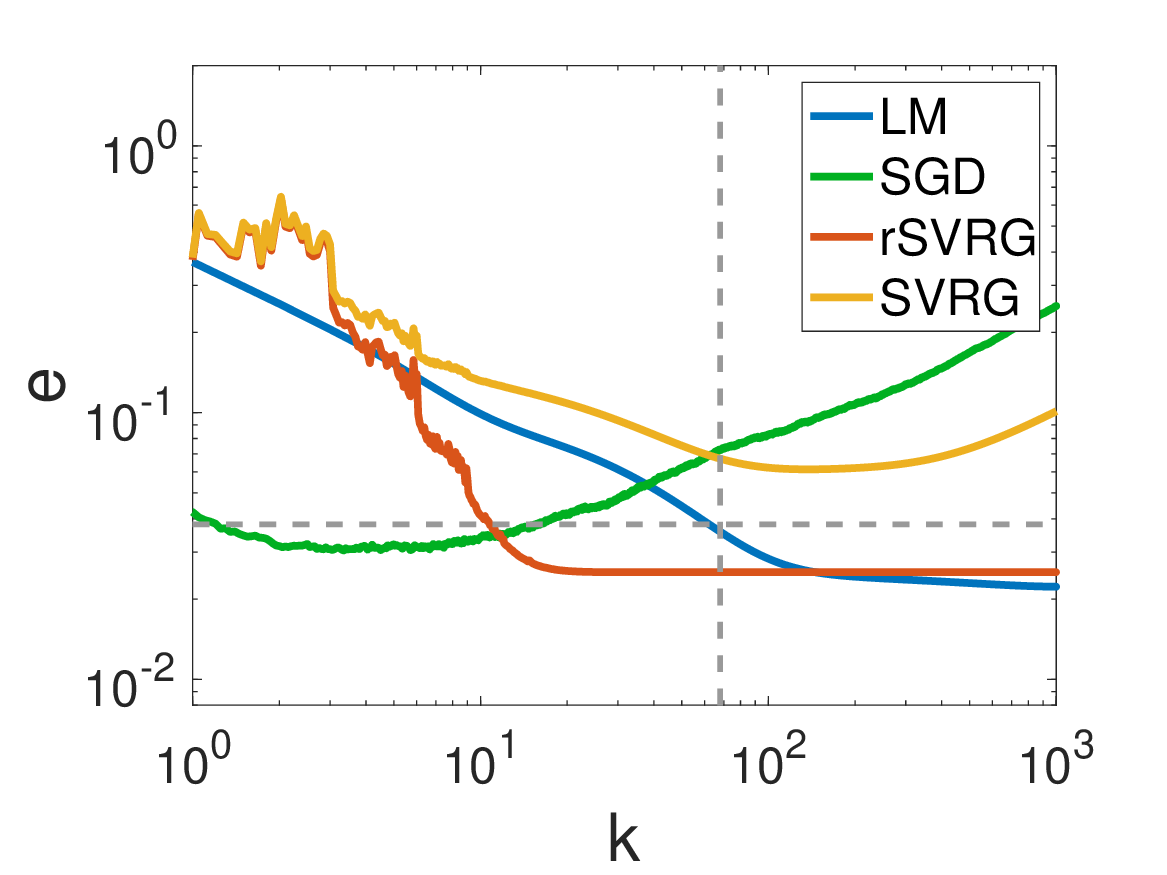}&
\includegraphics[width=0.31\textwidth,trim={1.5cm 0 0.5cm 0.5cm}]{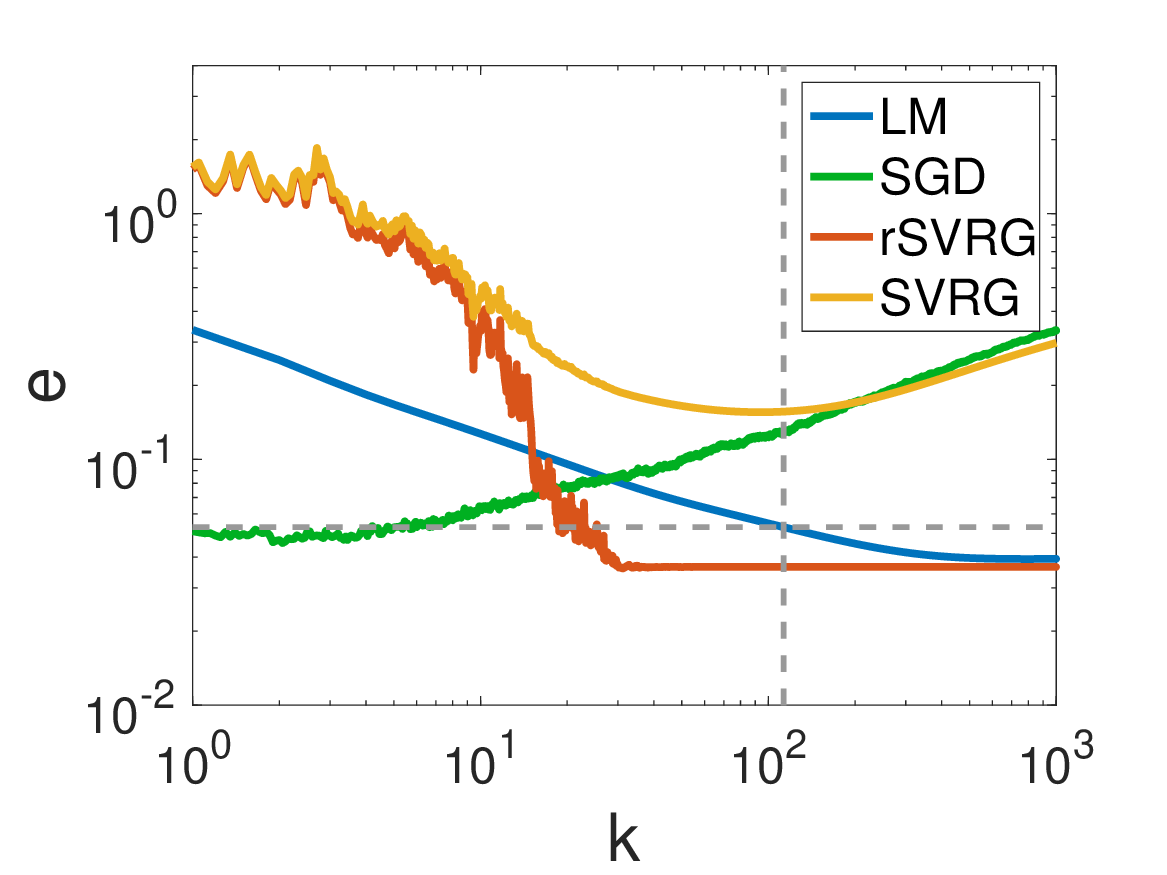}&
\includegraphics[width=0.31\textwidth,trim={1.5cm 0 0.5cm 0.5cm}]{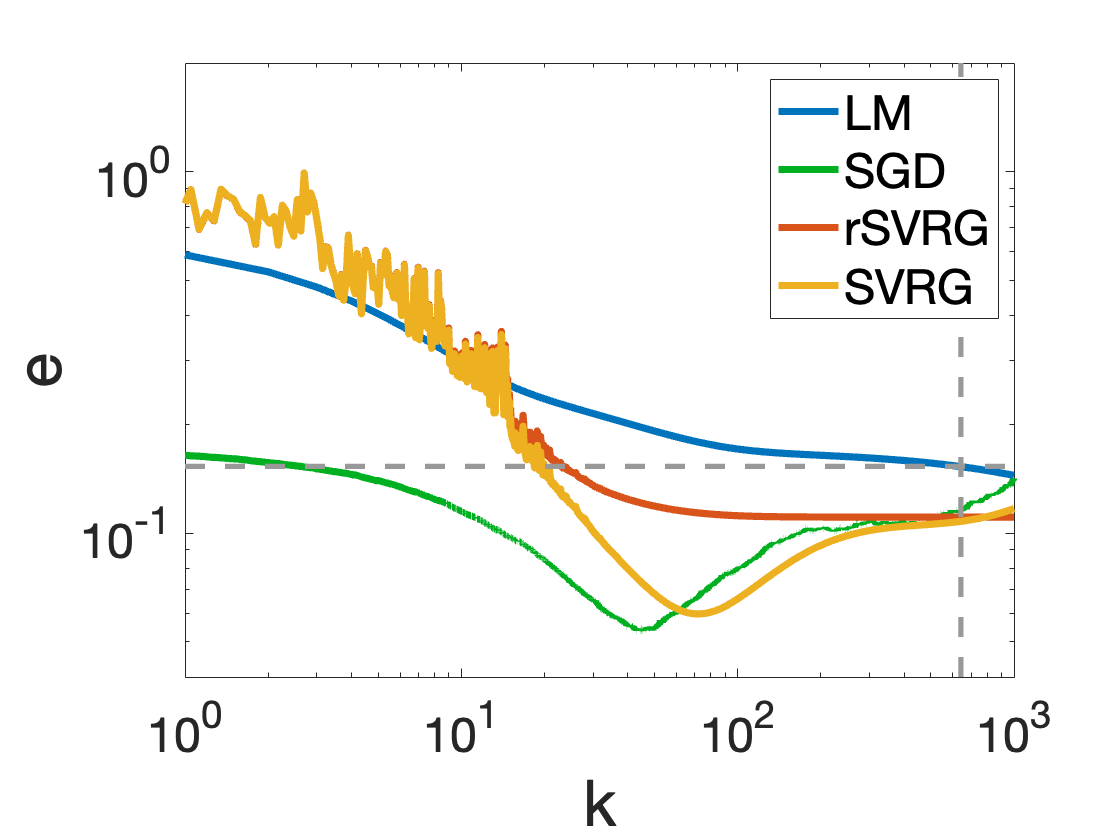}\\
\includegraphics[width=0.31\textwidth,trim={1.5cm 0 0.5cm 0.5cm}]{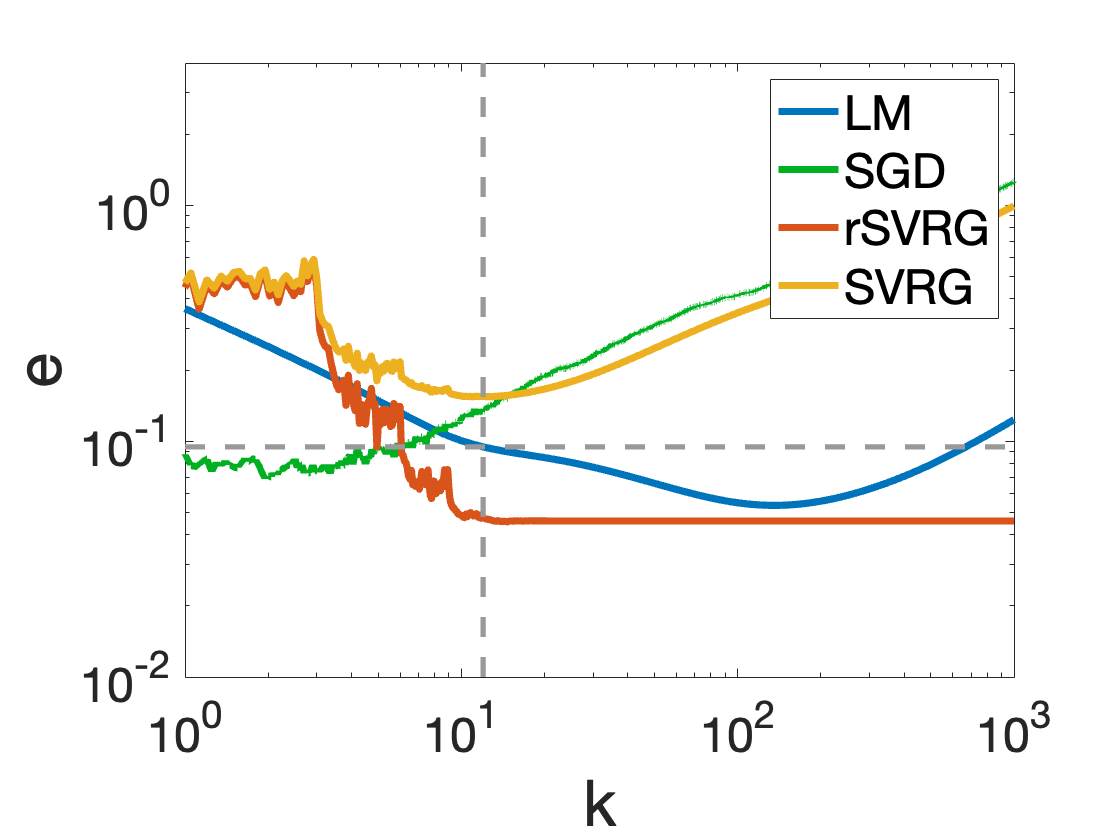}&
\includegraphics[width=0.31\textwidth,trim={1.5cm 0 0.5cm 0.5cm}]{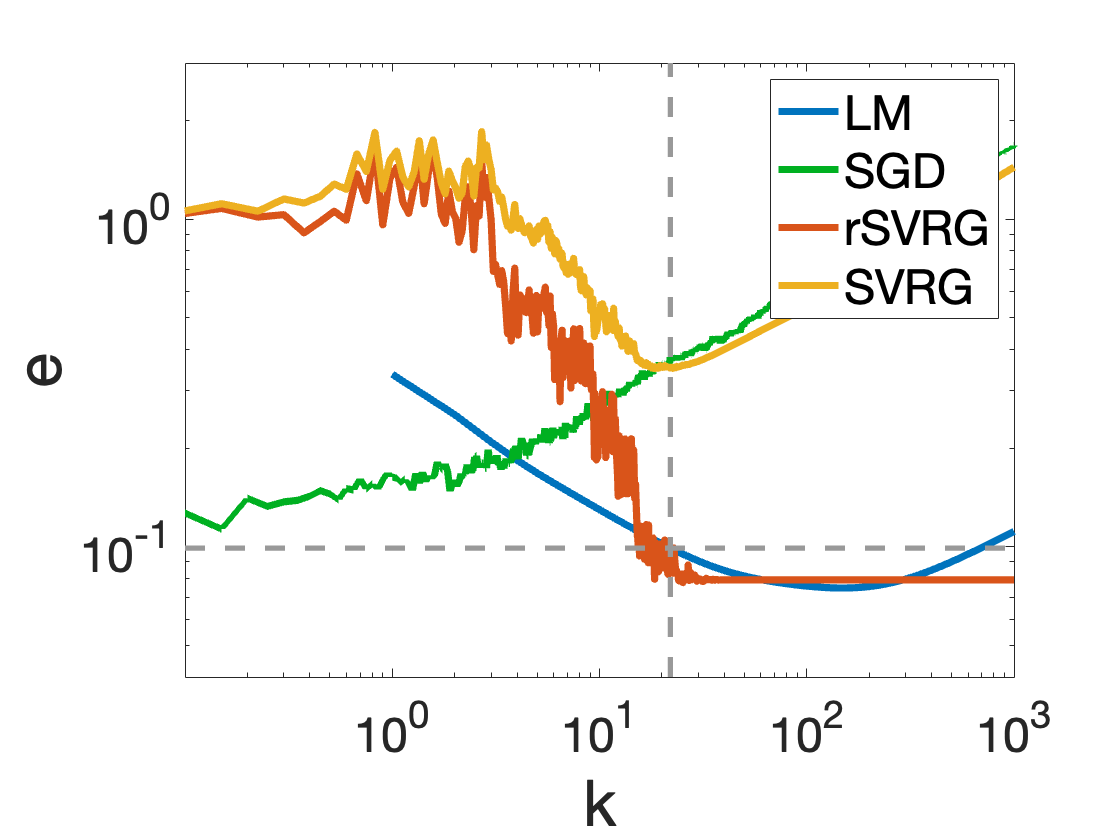}&
\includegraphics[width=0.31\textwidth,trim={1.5cm 0 0.5cm 0.5cm}]{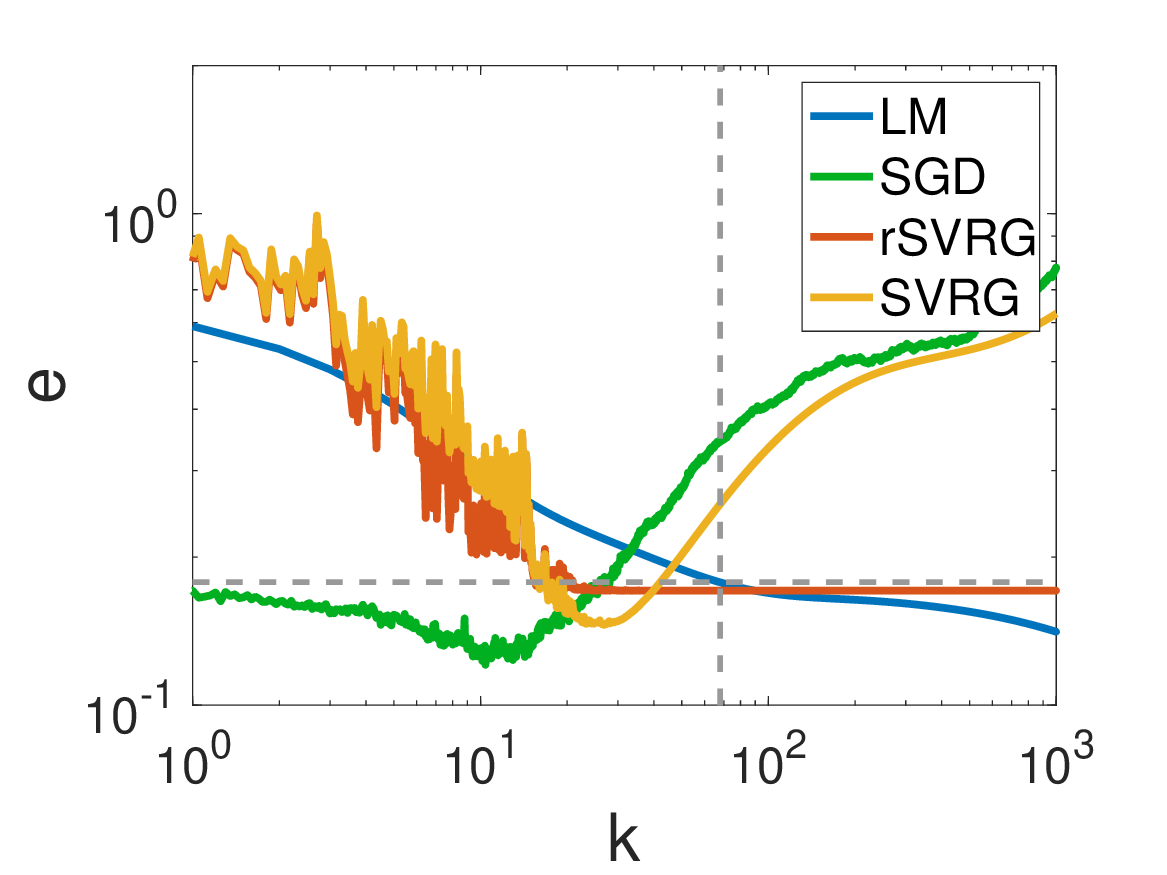}\\
\texttt{phillips}& \texttt{gravity} & \texttt{shaw}
\end{tabular}
\caption{The convergence of the relative error $e={\E[\| x_{k}^\delta-x_\dag\|^2]^\frac12}/{\|x_\dag\|}$ versus the iteration number $k$ for \texttt{phillips}, \texttt{gravity} and \texttt{shaw}. The rows from top to bottom are for $\epsilon=$1e-3,
$\epsilon=$5e-3, 
$\epsilon=$1e-2 and $\epsilon=$5e-2, respectively. The intersection of the gray dashed lines represents the stopping point, determined by the discrepancy principle, for LM along the iteration trajectory.}\label{fig}
\end{figure}

The numerical results for the  examples with varying regularity indices $\nu$ and noise levels $\epsilon$ are presented in Tables \ref{tab:phil}, \ref{tab:gravity}, and \ref{tab:shaw}.
It is observed that rSVRG achieves an accuracy (with much fewer iterations for relatively low-regularity problems) comparable to that for the LM across varying regularity. 
SVRG can also achieve comparable accuracy in low-regularity cases, indicating its optimality.
However, with the chosen step size schedule, it is not optimal for highly regular solutions, for which smaller step sizes are required to achieve the optimal  error \cite{JinZhouZou:2022ip}. While SGD can achieve optimal accuracy with the fewest iterations for problems with low-regularity solutions, its accuracy is no longer optimal for high-regularity cases (e.g., $\nu=1$). Moreover, (r)SVRG exhibits reduced oscillations near the optimum compared to SGD when the noise level $\epsilon$ is relatively large, cf. Fig. \ref{fig}. We refer interested readers to \cite{JinZhouZou:2022ip} for a more comprehensive comparison of the performance between SVRG and SGD under suitable step size schedules, as well as the numerical results that indicate the advantages of SVRG over LM.
Typically, for all the problems and the four methods, the case with a higher noise level requires fewer iterations to reach the desired accuracy.
These observations agree with the theoretical results of Theorem \ref{thm:main} and Corollary \ref{cor:rate}.
Moreover, the error of rSVRG at its plateau point is typically lower than that of the other three methods. Indeed, the results by rSVRG can be much more accurate (by a factor of 10) than that by SVRG for many settings, as indicated by the results in the tables. Furthermore, the required number of iterations is smaller for rSVRG than that for SVRG. These results clearly show the potential computational advantage of rSVRG (provided that the truncation level $J$ in rSVRG can be chosen properly).
The convergence trajectories of the methods for the examples with $\nu=0$ in Fig. \ref{fig} show the advantage of rSVRG over the other three benchmark methods as seen in Tables \ref{tab:phil}-\ref{tab:shaw}. 

\begin{table}[hbt!]
  \centering\small
  \begin{threeparttable}
  \caption{The comparison between (r)SVRG, SGD and LM for \texttt{s-phillips}.}\label{tab:phil}
    \begin{tabular}{ccccccccccccc}
    \toprule
    \multicolumn{2}{c}{Method}&
    \multicolumn{3}{c}{rSVRG  ($c_0=c/4$)}&\multicolumn{2}{c}{SVRG  ($c_0=c/4$)}& \multicolumn{2}{c}{SGD  ($c_0=c/4$)}&\multicolumn{2}{c}{LM}\\
    \cmidrule(lr){3-5} \cmidrule(lr){6-7} \cmidrule(lr){8-9}
    \cmidrule(lr){10-11}
    $\nu$& $\epsilon$ &$e_*$&$k_*$&$\lim_{k\to \infty} e$&$e_*$&$k_*$&$e_*$&$k_*$&$e_*$&$k_*$\\
    \hline
    $0$& 1e-3 & 1.93e-2 & 102.825 & 1.17e-2  & 1.52e-2 & 1170.900  & 9.48e-3 & 79.200  & 1.93e-2 & 758 \cr
       & 5e-3 & 2.81e-2 & 14.325  & 2.52e-2  & 6.13e-2 & 137.625   & 2.20e-2 & 13.400  & 2.81e-2 & 102 \cr
       & 1e-2 & 3.79e-2 & 12.000  & 2.63e-2  & 7.93e-2 & 70.050    & 3.04e-2 & 3.350   & 3.81e-2 & 68  \cr
       & 5e-2 & 8.81e-2 & 6.075   & 4.58e-2  & 1.54e-1 & 11.100    & 6.91e-2 & 1.950   & 9.44e-2 & 12  \cr
    \hline
    $0.25$& 1e-3 & 4.58e-3 & 206.700 & 4.29e-3  & 2.73e-2 & 819.225  & 6.22e-3 & 31.700  & 4.58e-3 & 135 \cr
          & 5e-3 & 1.48e-2 & 13.425  & 5.68e-3  & 5.73e-2 & 110.925  & 1.26e-2 & 4.950   & 1.48e-2 & 60  \cr
          & 1e-2 & 2.79e-2 & 12.825  & 9.43e-3  & 7.50e-2 & 58.650   & 1.78e-2 & 2.100   & 2.81e-2 & 26  \cr
          & 5e-2 & 4.13e-2 & 9.075   & 3.83e-2  & 1.37e-1 & 11.550   & 6.47e-2 & 1.750   & 4.66e-2 & 10  \cr
    \hline
    $0.5$& 1e-3 & 2.87e-3 & 24.300 & 1.01e-3  & 2.73e-2 & 841.575  & 6.45e-3 & 32.250  & 2.90e-3 & 94 \cr
         & 5e-3 & 1.00e-2 & 12.675 & 3.79e-3  & 5.79e-2 & 115.050  & 1.10e-2 & 3.250   & 1.21e-2 & 23 \cr
         & 1e-2 & 1.33e-2 & 11.475 & 7.52e-3  & 7.53e-2 & 60.375   & 1.61e-2 & 2.100   & 1.51e-2 & 16 \cr
         & 5e-2 & 2.85e-2 & 9.150  & 2.49e-2  & 1.44e-1 & 12.675   & 6.10e-2 & 0.600   & 2.92e-2 & 8  \cr
    \hline
    $1$& 1e-3 & 1.53e-3 & 15.225 & 7.22e-4  & 2.76e-2 & 866.250  & 5.79e-3 & 33.400  & 1.92e-3 & 25 \cr
       & 5e-3 & 3.35e-3 & 17.775 & 3.28e-3  & 5.93e-2 & 163.800  & 1.05e-2 & 3.250   & 3.44e-3 & 16 \cr
       & 1e-2 & 5.36e-3 & 14.700 & 4.36e-3  & 7.76e-2 & 66.900   & 1.52e-2 & 2.100   & 5.54e-3 & 12 \cr
       & 5e-2 & 1.57e-2 & 12.075 & 1.57e-2  & 1.43e-1 & 11.850   & 5.61e-2 & 0.350   & 1.82e-2 & 5  \cr
    \bottomrule
    \end{tabular}
    \end{threeparttable}
\end{table}

\begin{table}[htp!]
  \centering\small
  \begin{threeparttable}
  \caption{The comparison between (r)SVRG, SGD and LM for \texttt{s-gravity}.}\label{tab:gravity}
    \begin{tabular}{ccccccccccc}
    \toprule
    \multicolumn{2}{c}{Method}&
    \multicolumn{3}{c}{rSVRG ($c_0=c$)}&\multicolumn{2}{c}{SVRG ($c_0=c$)}& \multicolumn{2}{c}{SGD  ($c_0=c$)}&\multicolumn{2}{c}{LM}\\
    \cmidrule(lr){3-5} \cmidrule(lr){6-7} \cmidrule(lr){8-9}
    \cmidrule(lr){10-11}
    $\nu$& $\epsilon$ &$e_*$&$k_*$&$\lim_{k\to \infty} e$&$e_*$&$k_*$&$e_*$&$k_*$&$e_*$&$k_*$\\
    \hline
    $0$& 1e-3 & 2.36e-2 & 279.525 & 1.30e-2  & 4.12e-2 & 1356.150  & 1.66e-2 & 51.150  & 2.36e-2 & 1649 \cr
       & 5e-3 & 3.99e-2 & 32.325  & 2.33e-2  & 9.05e-2 & 247.650   & 3.49e-2 & 7.450   & 4.04e-2 & 255  \cr
       & 1e-2 & 4.93e-2 & 25.425  & 3.65e-2  & 1.56e-1 & 93.900    & 4.67e-2 & 2.050   & 5.30e-2 & 113  \cr
       & 5e-2 & 8.56e-2 & 22.950  & 7.92e-2  & 3.50e-1 & 18.450    & 1.13e-2 & 0.150   & 9.90e-2 & 22   \cr
    \hline
     $0.25$& 1e-3 & 6.16e-3 & 51.975 & 3.03e-3  & 4.74e-2 & 1550.400  & 7.52e-3 & 11.600  & 6.50e-3 & 319 \cr
           & 5e-3 & 1.56e-2 & 37.275 & 1.20e-2  & 1.25e-1 & 198.300   & 1.89e-2 & 2.100  & 1.64e-2 & 71  \cr
           & 1e-2 & 1.82e-2 & 27.150 & 1.27e-2  & 1.65e-1 & 164.325   & 3.21e-2 & 2.100  & 2.32e-2 & 43  \cr
           & 5e-2 & 5.12e-2 & 19.275 & 2.72e-2  & 4.05e-1 & 29.400    & 1.04e-1 & 0.150  & 5.35e-2 & 12  \cr
    \hline
    $0.5$& 1e-3 & 3.34e-3 & 44.625 & 2.31e-3  & 3.82e-2 & 1106.400  & 6.14e-3 & 11.350  & 3.39e-3 & 112 \cr
         & 5e-3 & 7.56e-3 & 47.025 & 5.52e-3  & 1.26e-1 & 206.325   & 1.77e-2 & 2.100   & 9.10e-3 & 40  \cr
         & 1e-2 & 1.33e-2 & 44.550 & 1.04e-2  & 1.59e-1 & 176.100   & 3.05e-2 & 2.100   & 1.41e-2 & 25  \cr
         & 5e-2 & 3.38e-2 & 20.925 & 1.02e-2  & 4.00e-1 & 29.400    & 1.20e-1 & 0.150   & 3.40e-2 & 8   \cr
    \hline
    $1$& 1e-3 & 1.41e-3 & 48.000 & 9.87e-4  & 3.82e-2 & 1222.725  & 6.00e-3 & 11.350  & 1.46e-3 & 42 \cr
       & 5e-3 & 3.06e-3 & 35.400 & 1.11e-3  & 1.07e-1 & 259.800   & 1.73e-2 & 2.100   & 4.11e-3 & 18 \cr
       & 1e-2 & 3.17e-3 & 33.000 & 1.43e-3  & 1.57e-1 & 161.175   & 2.97e-2 & 2.100   & 6.58e-3 & 12 \cr
       & 5e-2 & 1.08e-2 & 23.175 & 8.15e-3  & 3.92e-1 & 29.400    & 9.98e-2 & 0.150   & 1.48e-2 & 6  \cr
 \bottomrule
    \end{tabular}
    \end{threeparttable}
\end{table}

\begin{table}[htp!]
  \centering\small
  \begin{threeparttable}
  \caption{The comparison between (r)SVRG, SGD and LM for \texttt{s-shaw}.}\label{tab:shaw}
    \begin{tabular}{ccccccccccc}
    \toprule
    \multicolumn{2}{c}{Method}&
    \multicolumn{3}{c}{rSVRG ($c_0=c$)}&\multicolumn{2}{c}{SVRG ($c_0=c$)}& \multicolumn{2}{c}{SGD  ($c_0=c$)}&\multicolumn{2}{c}{LM}\\
    \cmidrule(lr){3-5} \cmidrule(lr){6-7} \cmidrule(lr){8-9}
    \cmidrule(lr){10-11}
    $\nu$& $\epsilon$ &$e_*$&$k_*$&$\lim_{k\to \infty} e$&$e_*$&$k_*$&$e_*$&$k_*$&$e_*$&$k_*$\\
    \hline
    $0$& 1e-3 & 4.94e-2 & 39.825 & 4.94e-2  & 3.41e-2 & 4183.950  & 3.41e-2 & 2732.000  & 4.93e-2 & 22314 \cr
       & 5e-3 & 9.22e-2 & 57.375 & 6.88e-2  & 4.93e-2 & 132.675   & 4.70e-2 & 69.400  & 9.28e-2 & 4858  \cr
       & 1e-2 & 1.53e-1 & 23.025 & 1.11e-1  & 5.98e-2 & 71.775    & 5.35e-2 & 45.400  & 1.53e-1 & 642   \cr
       & 5e-2 & 1.74e-1 & 20.925 & 1.71e-1  & 1.46e-1 & 26.925    & 1.21e-1 & 10.400  & 1.78e-1 & 68    \cr
    \hline
    $0.25$& 1e-3 & 1.69e-2 & 90.450 & 1.09e-2  & 2.01e-2 & 745.500  & 4.27e-3 & 45.200  & 1.69e-2 & 1218 \cr
          & 5e-3 & 2.21e-2 & 36.000 & 2.20e-2  & 4.34e-2 & 79.800   & 1.49e-2 & 11.500  & 2.24e-2 & 139  \cr
          & 1e-2 & 2.46e-2 & 23.625 & 2.24e-2  & 6.99e-2 & 56.550   & 2.26e-2 & 7.300   & 2.59e-2 & 99   \cr
          & 5e-2 & 5.21e-2 & 15.000 & 3.20e-2  & 1.75e-1 & 20.775   & 5.96e-2 & 0.900   & 7.02e-2 & 24   \cr
    \hline
    $0.5$& 1e-3 & 2.97e-3 & 42.075 & 2.84e-3  & 2.05e-2 & 598.725  & 4.73e-3 & 14.850  & 3.16e-3 & 169 \cr
         & 5e-3 & 7.80e-3 & 30.075 & 3.81e-3  & 5.17e-2 & 85.275   & 1.26e-2 & 4.150   & 8.83e-3 & 78  \cr
         & 1e-2 & 1.55e-2 & 21.075 & 5.89e-3  & 7.51e-2 & 56.175   & 1.72e-2 & 0.900   & 1.69e-2 & 42  \cr
         & 5e-2 & 4.63e-2 & 18.825 & 4.13e-2  & 1.97e-1 & 19.050   & 6.19e-2 & 0.900   & 5.36e-2 & 16  \cr
    \hline
    $1$& 1e-3 & 1.60e-3 & 40.875 & 5.63e-4  & 2.07e-2 & 225.300  & 5.31e-3 & 13.800  & 1.80e-3 & 54 \cr
       & 5e-3 & 5.16e-3 & 41.475 & 2.81e-3  & 5.60e-2 & 84.075   & 1.28e-2 & 0.900   & 6.13e-3 & 25 \cr
       & 1e-2 & 7.24e-3 & 28.650 & 6.31e-3  & 8.20e-2 & 55.125   & 1.75e-2 & 0.900   & 1.18e-2 & 19 \cr
       & 5e-2 & 4.79e-2 & 18.000 & 1.91e-2  & 2.12e-1 & 16.650   & 6.68e-2 & 0.900   & 5.26e-2 &  6 \cr
 \bottomrule
    \end{tabular}
    \end{threeparttable}
\end{table}

\begin{figure}[hbt!]
\centering
  \setlength{\tabcolsep}{4pt}
\begin{tabular}{ccc}
\includegraphics[width=0.31\textwidth,trim={1.5cm 0 0.5cm 0.5cm}]{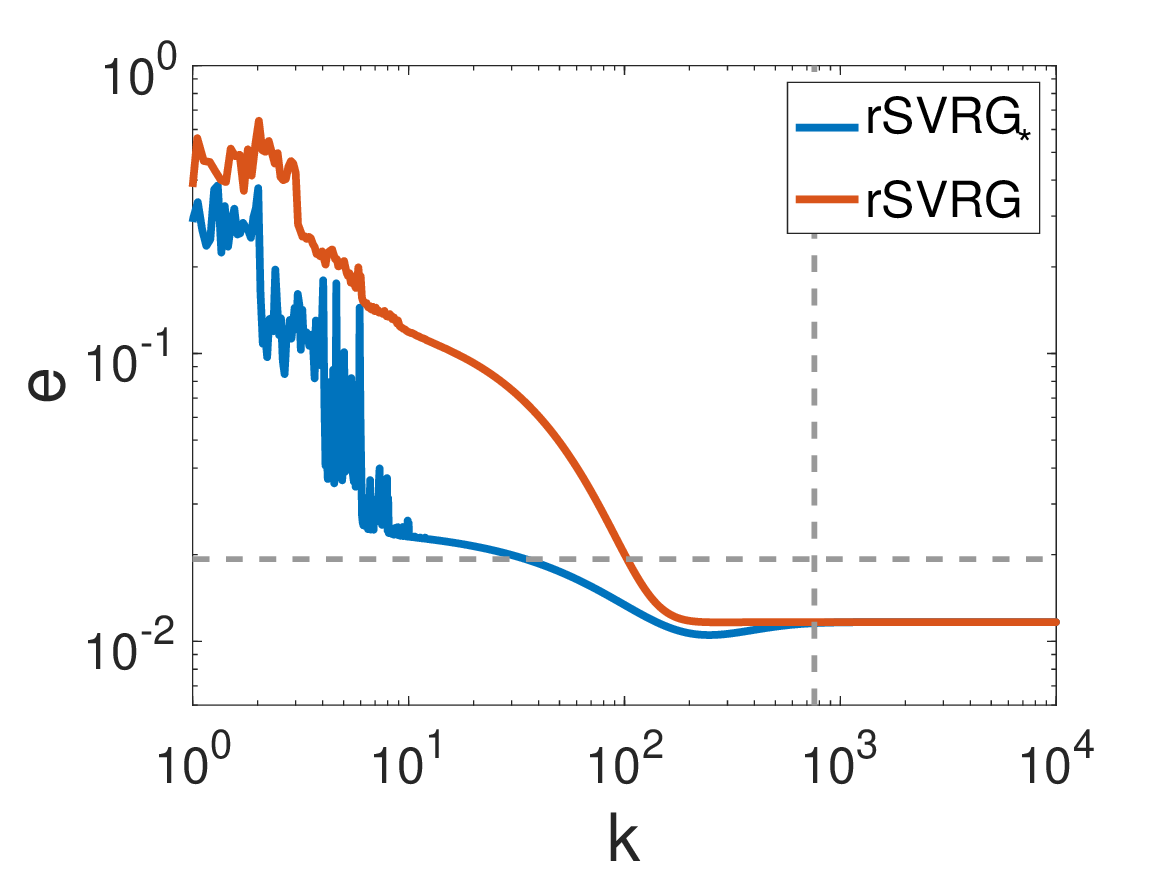}&
\includegraphics[width=0.31\textwidth,trim={1.5cm 0 0.5cm 0.5cm}]{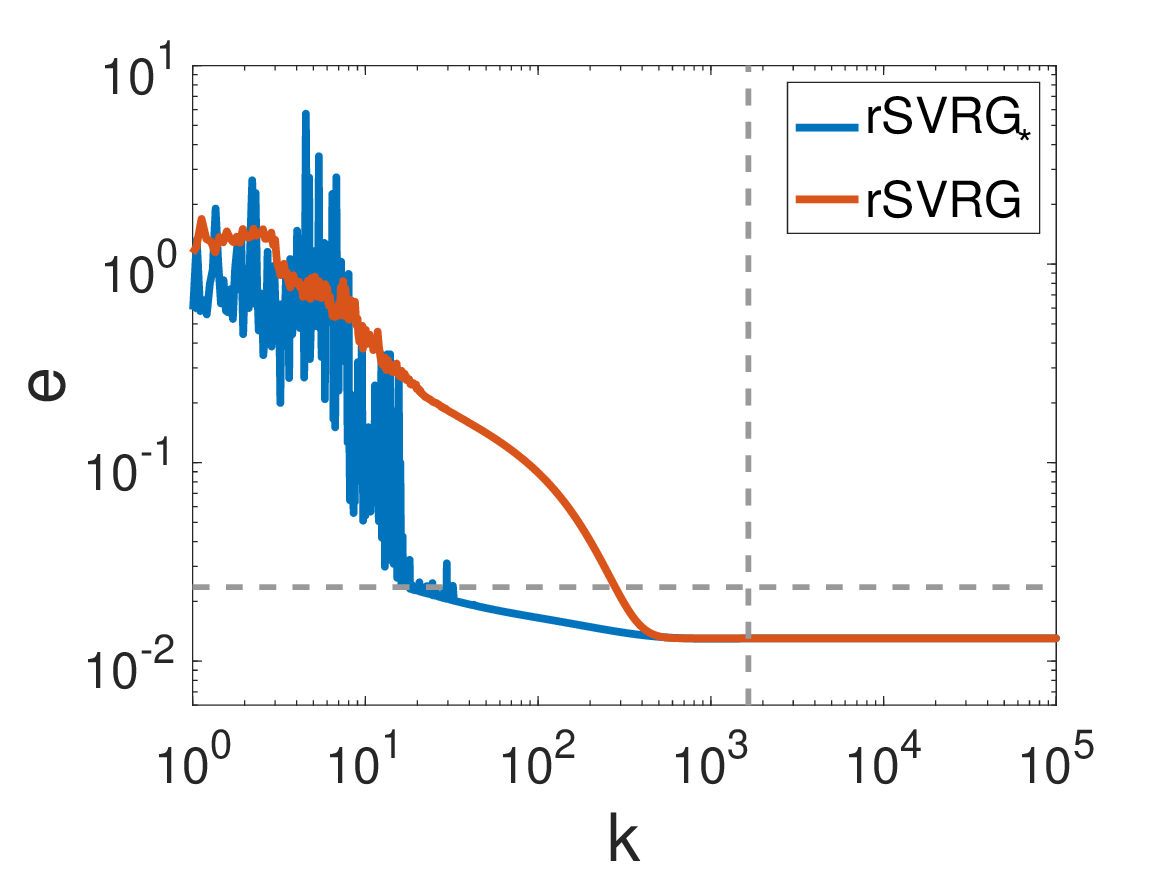}&
\includegraphics[width=0.31\textwidth,trim={1.5cm 0 0.5cm 0.5cm}]{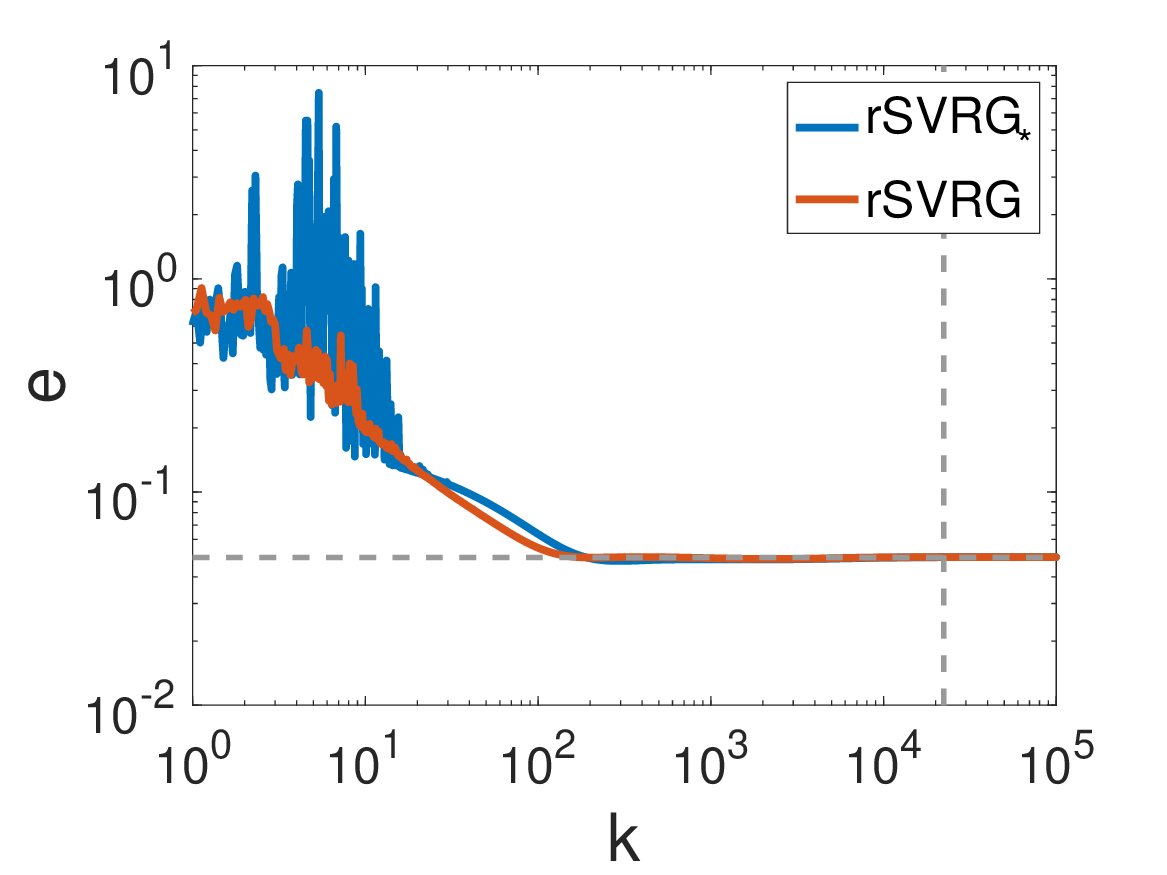}\\
\includegraphics[width=0.31\textwidth,trim={1.5cm 0 0.5cm 0.5cm}]{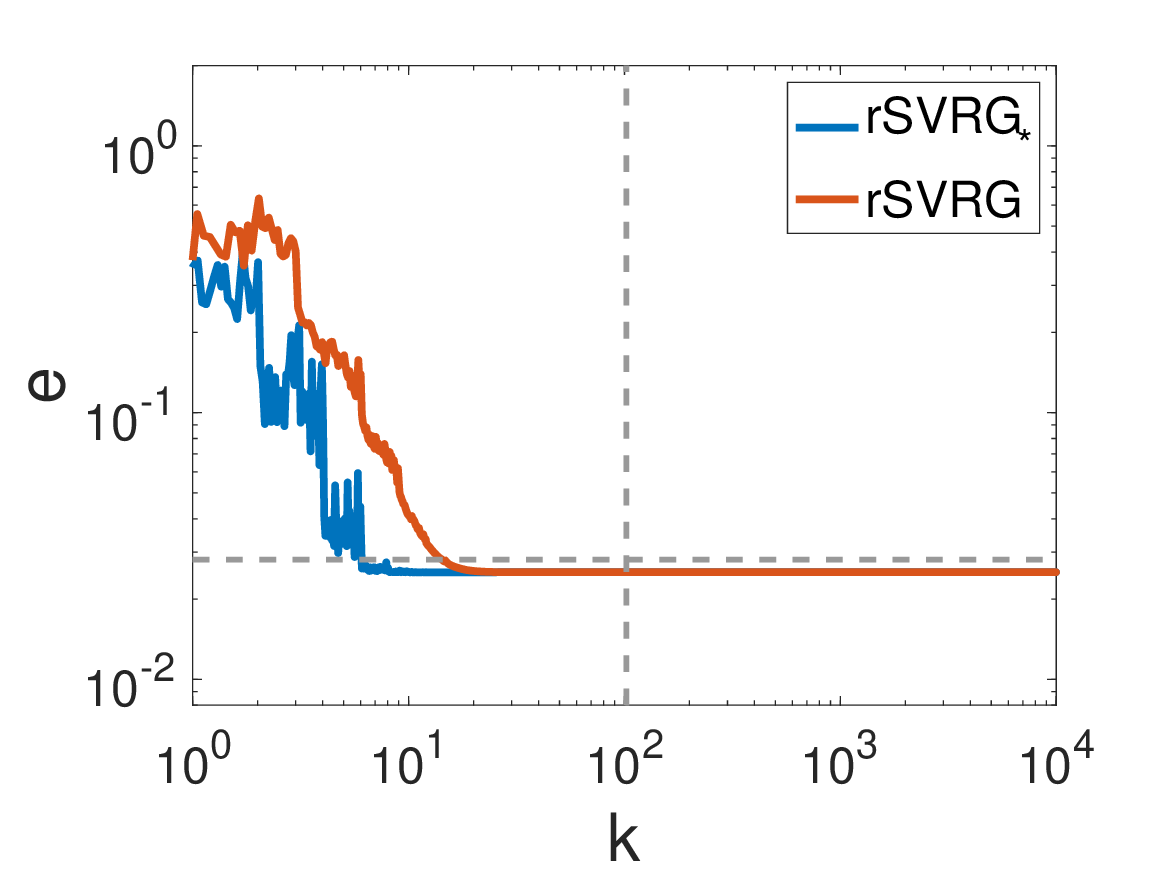}&
\includegraphics[width=0.31\textwidth,trim={1.5cm 0 0.5cm 0.5cm}]{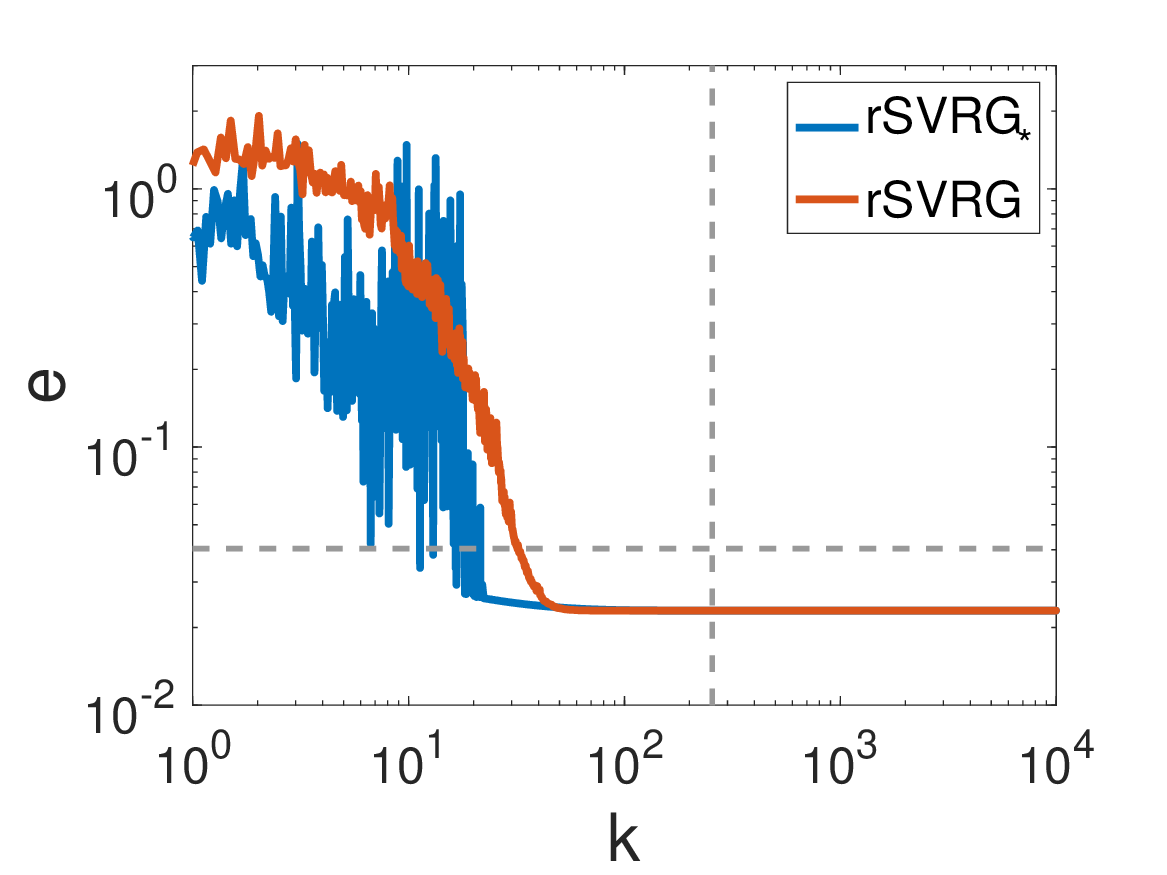}&
\includegraphics[width=0.31\textwidth,trim={1.5cm 0 0.5cm 0.5cm}]{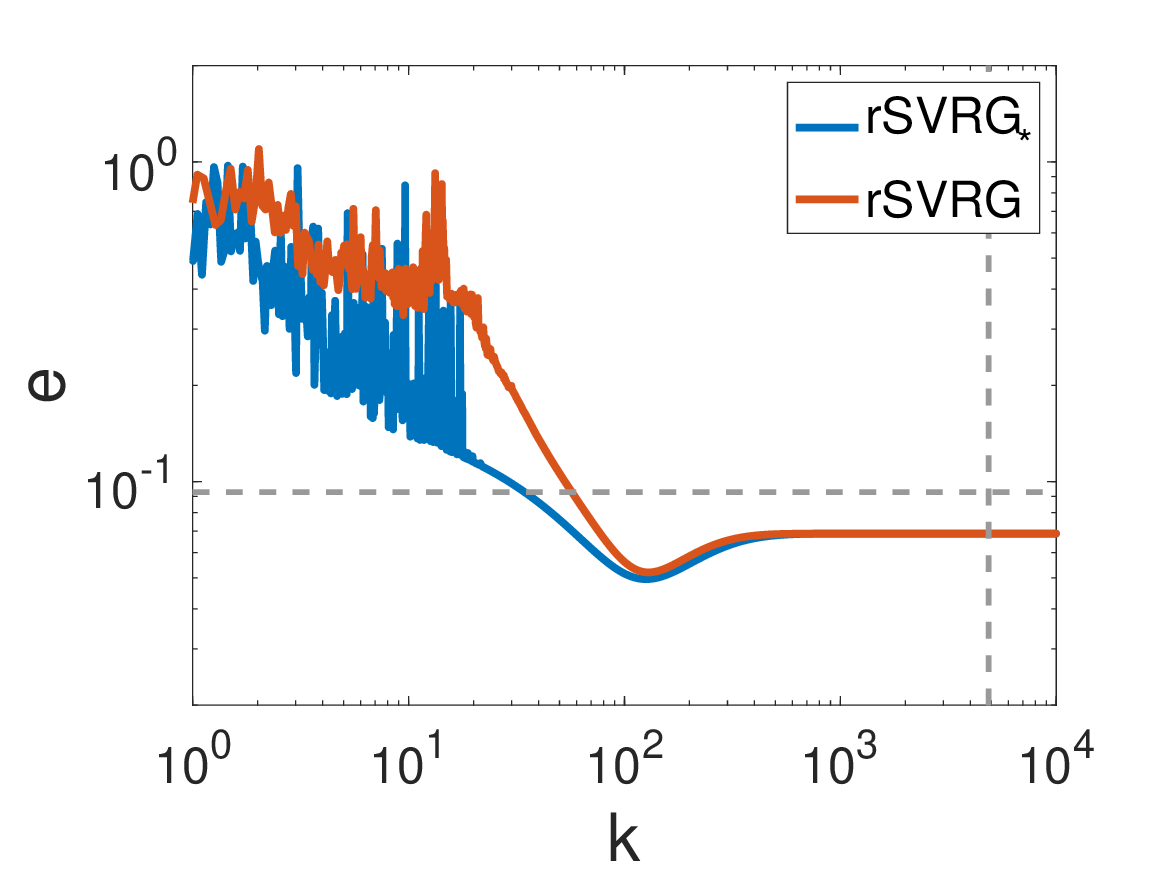}\\
\includegraphics[width=0.31\textwidth,trim={1.5cm 0 0.5cm 0.5cm}]{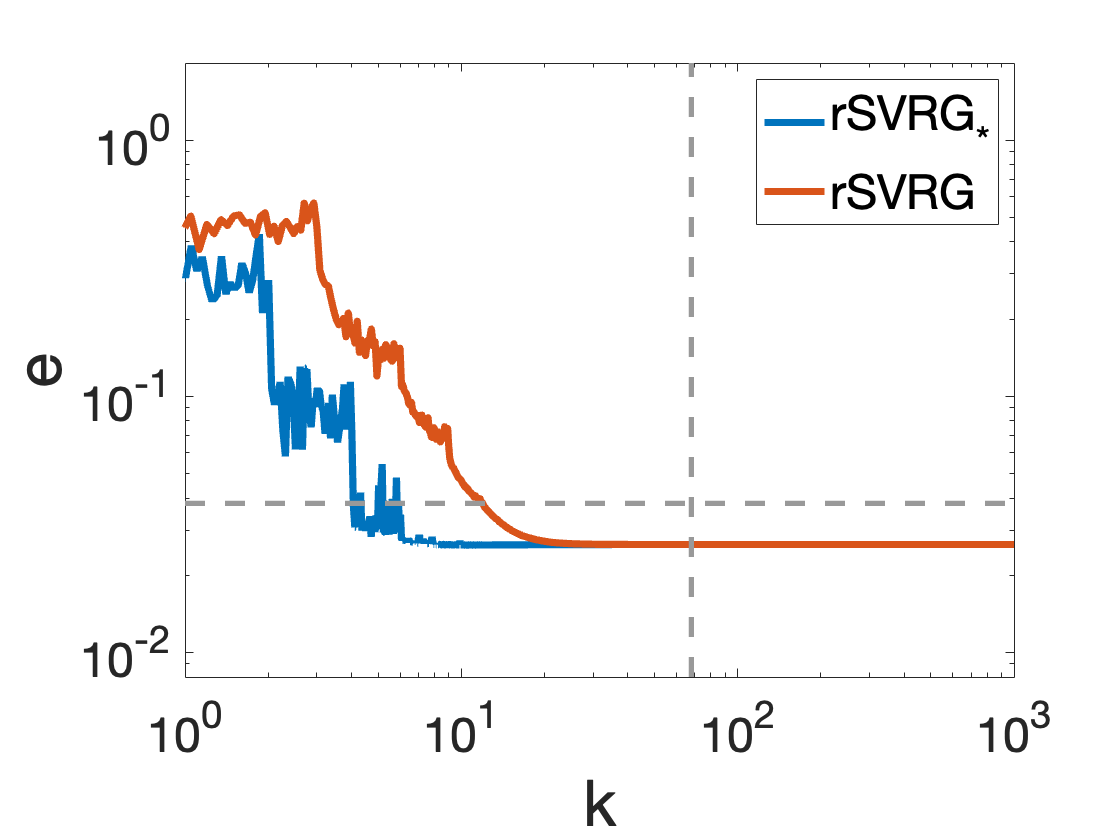}&
\includegraphics[width=0.31\textwidth,trim={1.5cm 0 0.5cm 0.5cm}]{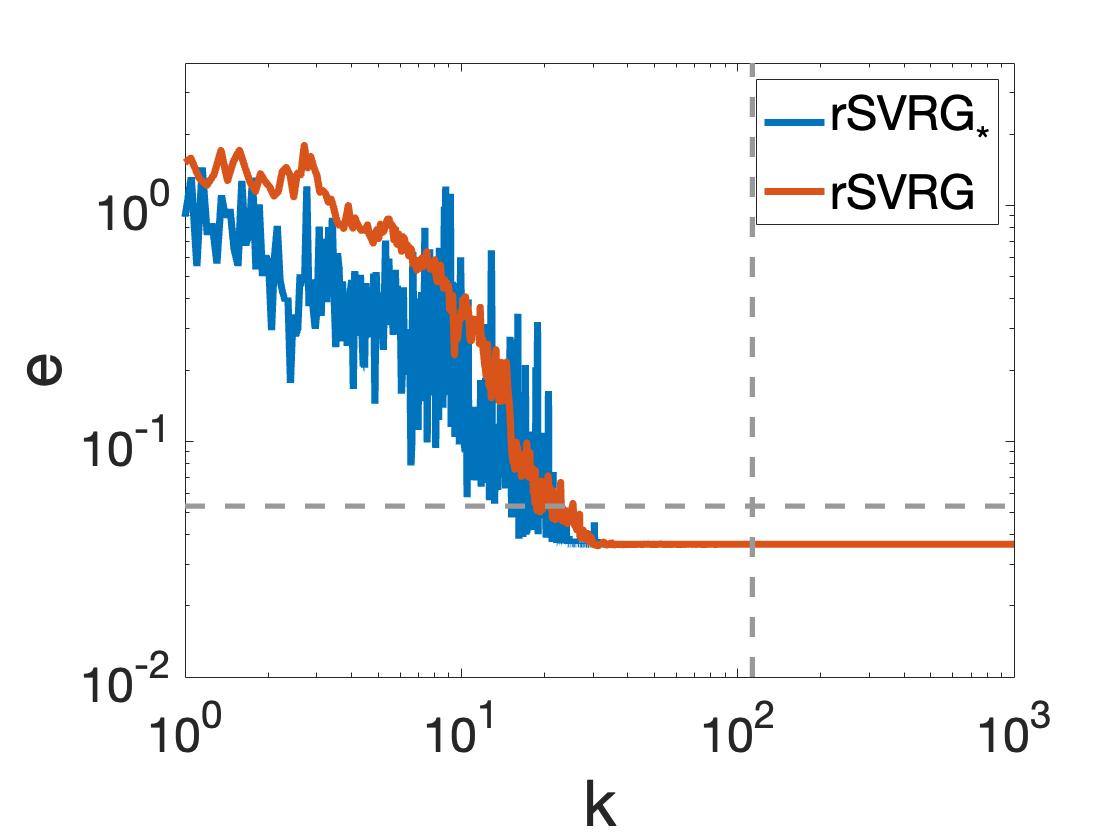}&
\includegraphics[width=0.31\textwidth,trim={1.5cm 0 0.5cm 0.5cm}]{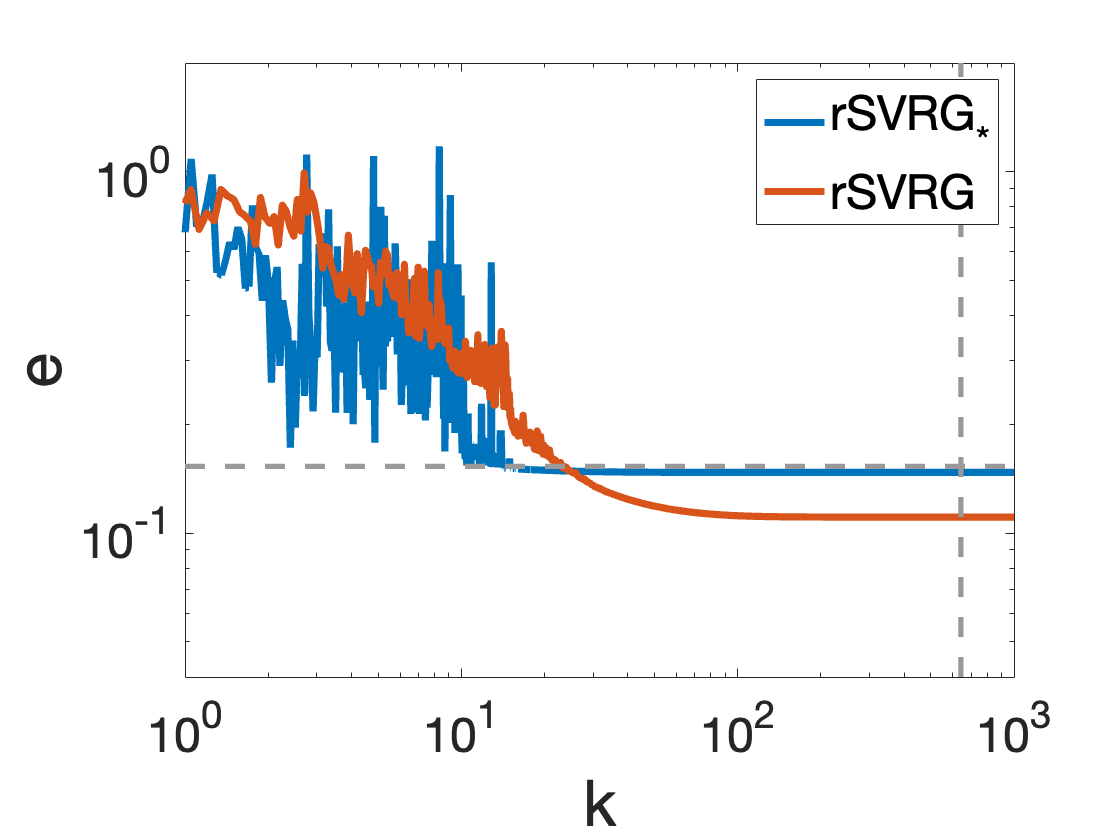}\\
\includegraphics[width=0.31\textwidth,trim={1.5cm 0 0.5cm 0.5cm}]{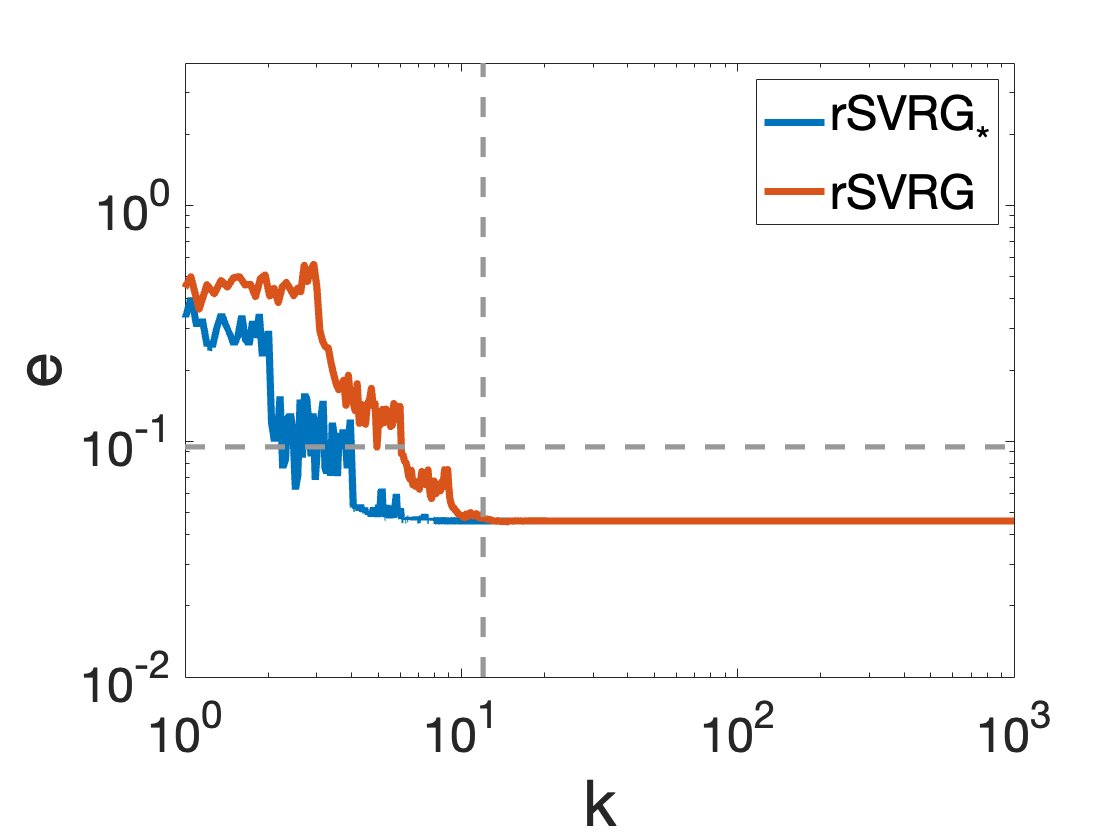}&
\includegraphics[width=0.31\textwidth,trim={1.5cm 0 0.5cm 0.5cm}]{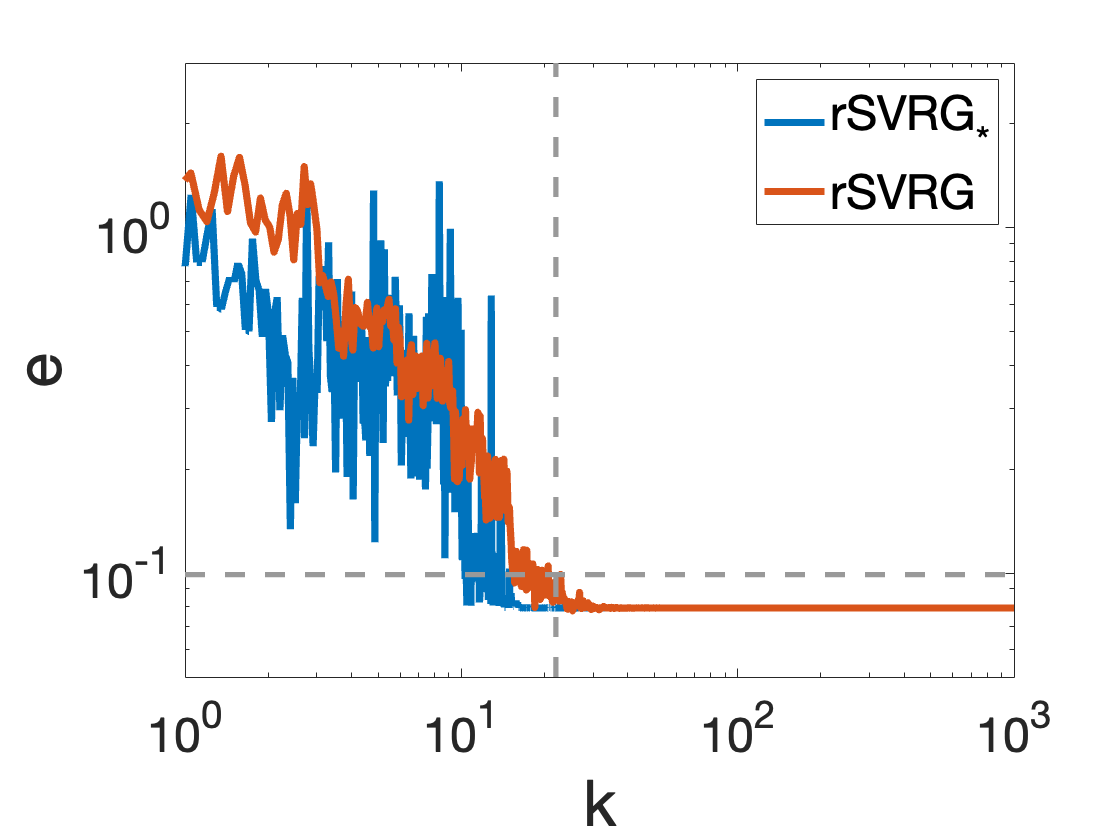}&
\includegraphics[width=0.31\textwidth,trim={1.5cm 0 0.5cm 0.5cm}]{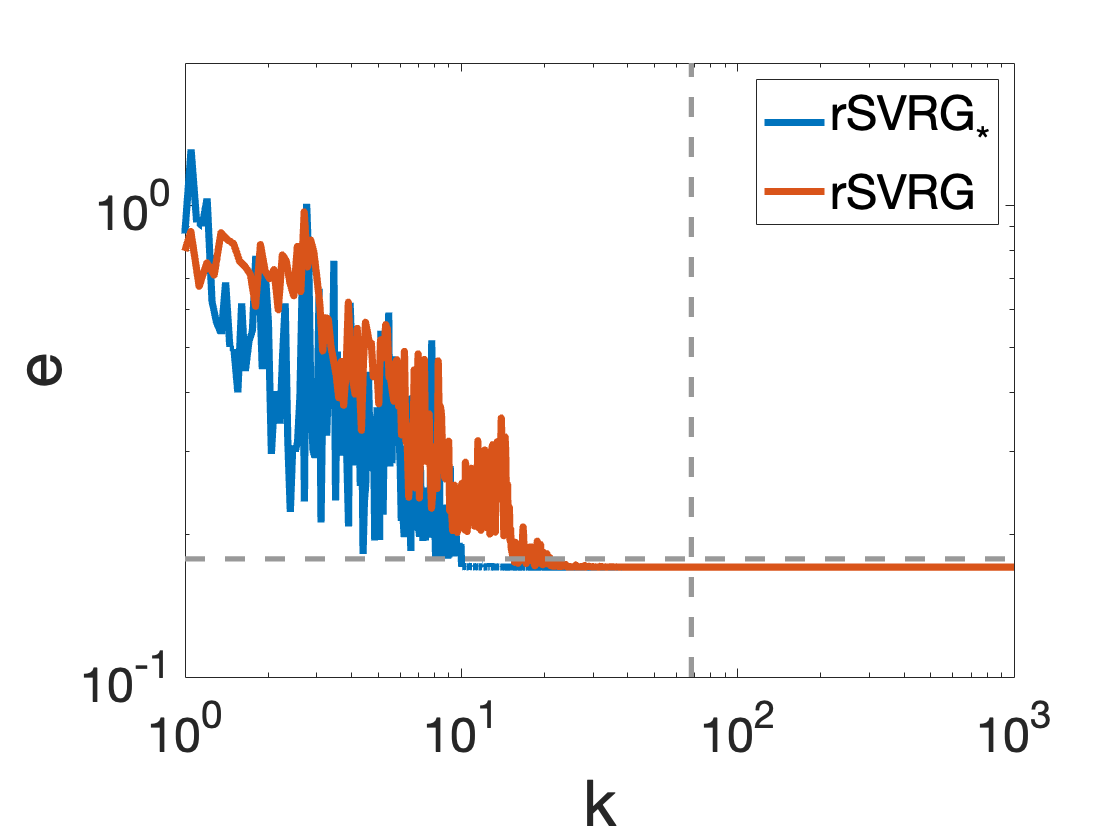}\\
\texttt{phillips}& \texttt{gravity} & \texttt{shaw}
\end{tabular}
\caption{The convergence of the relative error $e={\E[\| x_{k}^\delta-x_\dag\|^2]^\frac12}/{\|x_\dag\|}$ versus the iteration number $k$ for \texttt{phillips}, \texttt{gravity} and \texttt{shaw}. The rows from top to bottom are for $\epsilon=$1e-3,
$\epsilon=$5e-3, 
$\epsilon=$1e-2 and $\epsilon=$5e-2, respectively. The intersection of the gray dashed lines represents the stopping point, determined by the discrepancy principle, for LM along the iteration trajectory.}\label{fig-svrg}
\end{figure}

In addition, to shed further insights into rSVRG, we have also implemented an alternative variant, denoted by 
$\mbox{rSVRG}_*$, which uses the SVD of $A$, as described in the introduction. The constant step size is taken to be $c_0=\mathcal{O}(c_*)$ for $\mbox{rSVRG}_*$, with $c_*=\min_{1\leq i\leq J}(\sigma_i^{-2})$ ($J$ denotes the truncation level, 
and is much smaller than the data size $n$). The relevant numerical results are presented in Fig. \ref{fig-svrg}, in which we compare the rSVRG in Algorithm \ref{alg:dsvrg} with rSVRG$_*$, where one epoch refers to $nM/(J+M)$ rSVRG$_*$ iterations. Across all the four noise levels, we consistently observe that the iteration trajectories of the variant rSVRG$_*$ exhibit notable oscillations during the initial phase, which can be rather pronounced. However, the oscillations nearly disappear when the iterates are getting closer to the convergence. Moreover, for relatively mildly ill-posed problems, rSVRG$_*$ reaches the convergence slightly faster than rSVRG, but the final converged values are visually indistinguishable. Thus, the two variants (rSVRG and rSVRG$_*$) perform more or less comparably with each other.

\section{Concluding remarks} \label{sec:conc}

In this work, we have investigated stochastic variance reduced gradient (SVRG) and a regularized variant (rSVRG) for solving linear inverse problems in Hilbert spaces. 
We have established the regularizing property of both SVRG and rSVRG. Under the source condition, we have derived convergence rates in expectation and in the uniform sense for (r)SVRG. These results indicate the optimality of SVRG for nonsmooth solutions and the built-in regularization mechanism and optimality of rSVRG. The numerical results for three linear inverse problems with varying degree of ill-posedness show the advantages of rSVRG over SGD, standard SVRG and Landweber method. Note that both SVRG and rSVRG depend on the knowledge of the noise level, and the \textit{a priori} choice of the relevant algorithmic parameters also require a knowledge of the unknown regularity index $\nu$. Thus there is a need to develop data-driven or adaptive rules for choosing these algorithmic parameters. Moreover, in practice, the noise level may be unknown, and certain heuristic techniques are required for their efficient implementation, e.g., as the \textit{a priori} stopping rule or constructing the approximate operator $A$. 
We leave these interesting questions to future works.

\section*{Acknowledgments} 
The authors are grateful to the editor and two anonymous referees for their constructive comments which have improved the quality of the paper.

\appendix

\section{Proof of Theorem \ref{thm:Delta}}\label{app:estimate}
In this part, we give the technical proof of Theorem \ref{thm:Delta}. First we give two technical estimates.

\begin{lemma}\label{lem:kernel}
Under Assumption \ref{ass}(i), for any $s\geq 0$, $k,t\in \mathbb{N}$ and $\epsilon\in(0,\frac12]$, there hold
\begin{align*}
&\|B^s P^k\|\leq s^s c_0^{-s}  k^{-s}, \quad
\|(I-P^t)P^k\|\leq t (k+t)^{-1},\\
& \|B^\frac12(I-P^t)P^k\|\leq 2^{1+\epsilon}\epsilon^\epsilon c_0^{-\epsilon}\|B\|^{\frac12-\epsilon}  t 
(k+t)^{-(1+\epsilon)}.
\end{align*}
\end{lemma}
\begin{proof}
The first inequality can be found in \cite[Lemma 3.4]{JinZhouZou:2022ip}. To show the second inequality, let $\mathrm{Sp}(P)$ be the spectrum of $P$. Then there holds
\begin{align*}
\|(I-P^t)P^k\|=&\sup_{\lambda\in{\rm Sp}(P)}|(1-\lambda^t)\lambda^k|\leq \sup_{\lambda\in[0,1]}(1-\lambda^t)\lambda^k.
\end{align*}
Let $g(\lambda)=(1-\lambda^t)\lambda^k$. Then $g'(\lambda)=k\lambda^{k-1}-(k+t)\lambda^{k+t-1}$, so that
$g(\lambda)$ achieves its maximum over the interval $[0,1]$ at $\lambda=\lambda_*$ with $\lambda_*^t=\frac{k}{k+t}=1-\frac{t}{k+t}$. Consequently,
\begin{align*}
\|(I-P^t)P^k\|\leq g(\lambda_*)
\leq t (k+t)^{-1}.
\end{align*}
The last one follows by
\begin{align*}
\|B^\frac12(I-P^t)P^k\|&\leq \|B\|^{\frac12-\epsilon}\|B^\epsilon P^\frac{k+t}{2}\|\|(I-P^t)P^\frac{k-t}{2}\|\\
&\leq 2^{1+\epsilon}\epsilon^\epsilon c_0^{-\epsilon}\|B\|^{\frac12-\epsilon}  t 
(k+t)^{-(1+\epsilon)}.
\end{align*}
This completes the proof of the lemma.
\end{proof}

\begin{lemma}\label{lem:N}
Let $R: X\rightarrow X$ be a deterministic bounded linear operator. Then for any $j\geq 0$, there hold
\begin{align*}
\E[\|R N_j \Delta_j^\delta\|^2|\mathcal{F}_j]\leq& \min\big(n^{-1}L\|R\|^2 ,\|R B^\frac12\|^2 \big)\|A\Delta_j^\delta\|^2,\\
\|R N_j \Delta_j^\delta\|\leq& \min\big(\sqrt{L}\|R\| ,\sqrt{n}\|R B^\frac12\| \big)\|A\Delta_j^\delta\|.
\end{align*}
\end{lemma}
\begin{proof}
The definitions of $N_j$ and $ B=\E[A_{i_j}^*A_{i_j}|\mathcal{F}_j]$ and the bias-variance decomposition imply
\begin{align*}
&\E[\|R N_j \Delta_j^\delta\|^2|\mathcal{F}_j]=\E[\|R (B-A_{i_j}^*A_{i_j}) \Delta_j^\delta\|^2|\mathcal{F}_j]\\
=&\E[\|R A_{i_j}^*A_{i_j} \Delta_j^\delta\|^2|\mathcal{F}_j]-\|R B \Delta_j^\delta\|^2
\leq L\|R\|^2  \E[\|A_{i_j} \Delta_j^\delta\|^2|\mathcal{F}_j]\\
=& L\|R\|^2n^{-1} \sum_{i=1}^n\|A_{i} \Delta_j^\delta\|^2=n^{-1} L\|R\|^2 \|A \Delta_j^\delta\|^2.
\end{align*}
Note that $N_j=A^*(n^{-1}A-b_{i_j}A_{i_j})$, with $b_{i_j}\in\mathbb{R}^n$ being the ${i_j}$th Cartesian basis vector.
Then the identity $\E[b_{i_j}A_{i_j} \Delta_j^\delta|\mathcal{F}_j]=n^{-1}A\Delta_j^\delta$ and the bias-variance decomposition yield
\begin{align*}
&\E[\|R N_j \Delta_j^\delta\|^2|\mathcal{F}_j]=\E[\|R A^*(n^{-1}A-b_{i_j}A_{i_j}) \Delta_j^\delta\|^2|\mathcal{F}_j]\\
\leq &\E[\|R A^*\|^2\|(n^{-1}A-b_{i_j}A_{i_j}) \Delta_j^\delta\|^2|\mathcal{F}_j]\\
=&n\|R B^\frac12\|^2\big(\E[\|b_{i_j}A_{i_j}\Delta_j^\delta\|^2|\mathcal{F}_j]-\|n^{-1}A \Delta_j^\delta\|^2\big)\\ \leq &n\|R B^\frac12\|^2\E[\|A_{i_j}\Delta_j^\delta\|^2|\mathcal{F}_j]
\leq
\|R B^\frac12\|^2 \|A\Delta_j^\delta\|^2.
\end{align*}
These estimates and the inequality 
$$\|R N_j \Delta_j^\delta\|^2\leq n \E[\|R N_j \Delta_j^\delta\|^2|\mathcal{F}_j] $$ complete the proof of the lemma.
\end{proof}

The proof of Theorem \ref{thm:Delta} is lengthy and very technical, and requires several technical lemmas with sharp constants. The first lemma provides useful bounds on the bias and variance components of the weighted successive error $A\Delta_{k}^\delta$ in terms of the iteration index $k$. 

\begin{lemma}\label{lem:decom_Delta}
Let Assumption \ref{ass}(i) hold. Then for any $k\geq 1$, $k_{\rm c}:=k-M-2$, $k_M:=[k/M]M$ and $\epsilon\in(0,\frac12]$, there hold
\begin{align}
\|\E[A\Delta_{k}^\delta]\|&\leq (\|A\|\|e_0^\delta\|+ \delta) M k^{-1},\label{eqn:Delta_E1}\\
\E[\|A\Delta_{k}^\delta-\E[A\Delta_{k}^\delta]\|^2]&\leq 
\left\{
\begin{aligned}
n^{-1}c_0^2 L\|A\|^2\Big(M^2\overline{\Phi}_{1}^{k_{\rm c}}(M+1,2)+\overline{\Phi}_{k_{\rm c}+1}^{k-1}(0,0)\Big),\\
2^{-1}c_0 L\Big(8M^2\overline{\Phi}_{1}^{k_{\rm c}}(M+1,3)+\overline{\Phi}_{k_{\rm c}+1}^{k_M-1}(0,1)+\overline{\Phi}_{k_M+1}^{k-1}(0,1)\Big),
\end{aligned}
\right.\label{eqn:Delta_E2}\\
\|A\Delta_{k}^\delta-\E[A\Delta_{k}^\delta]\|&\leq  c_0\sqrt{L}\|A\|\Big(\tfrac{2^{1+\epsilon}\epsilon^\epsilon n^{\epsilon}}{c_0^{\epsilon}\|A\|^{2\epsilon}}M \Phi_{1}^{k_{\rm c}}(M+1,1+\epsilon)+\Phi_{k_{\rm c}+1}^{k-1}(0,0)\Big).\label{eqn:Delta_E3}
\end{align}
\end{lemma}
\begin{proof}
Let $k=KM+t$ with $K\geq0$ and $1\leq t\leq M-1$. Similar to the proof of Lemma \ref{lem:decom_err}, 
for the bias term $\|\E[A\Delta_{KM+t}^\delta]\|$, by the definitions of $\zeta$, $B$ and $\tilde\phi^{t-1}$, the identity $$\|AP^{KM}\tilde\phi^{t-1} A^*\|=nc_0^{-1}\|(I-P^t)P^{KM}\|,$$ 
and Lemma \ref{lem:kernel}, we derive the estimate \eqref{eqn:Delta_E1} from Lemma \ref{lem:bias-var} that
\begin{align*}
&\|\E[A\Delta_{KM+t}^\delta]\|\leq \|A(P^t-I)P^{KM}\|\|e_0^\delta\|+n^{-1}c_0 \|AP^{KM}\tilde\phi^{t-1} A^*\| \delta\nonumber\\
\leq &(\|A\|\|e_0^\delta\|+ \delta)\|(I-P^t)P^{KM}\|\leq (\|A\|\|e_0^\delta\|+ \delta)M(KM+t)^{-1}.
\end{align*}
Next let 
\begin{equation*} 
S_{t,j}=c_0\|A(P^t-I)P^{KM-1-j} N_j\Delta_j^\delta\|\quad \mbox{and}\quad T_{t,j}=c_0\|AP^{KM+t-1-j} N_j\Delta_j^\delta\|.
\end{equation*}
Then for the variance, when $K\geq 1$, by Lemma \ref{lem:bias-var} and the identity $\E[\langle R_1 N_i \Delta_i^\delta, R_2 N_j \Delta_j^\delta\rangle|\mathcal{F}_j]=0$ for any $j>i$ and deterministic bounded linear operators $R_1$ and $R_2$, we have 
\begin{align*}
&\E[\|A\Delta_{KM+t}^\delta-\E[A\Delta_{KM+t}^\delta]\|^2]={\rm I}_1+{\rm I}_2+{\rm I}_3,
\end{align*}
with the three terms ${\rm I}_1$, ${\rm I}_2$ and ${\rm I}_3$ given by
\begin{align*}
 {\rm I}_1=\sum_{j=1}^{(K-1)M+t-2}\E [S_{t,j}^2], \quad 
{\rm I}_2=\sum_{j=(K-1)M+t-1}^{KM-1} \E[S_{t,j}^2]\quad \mbox{and}\quad
{\rm I}_3=\sum_{j=KM+1}^{KM+t-1}\E[T_{t,j}^2].  
\end{align*}
By Lemma \ref{lem:N}, the following estimates hold
\begin{align*}
    \E[S_{t,j}^2]&\leq c_0^2L\|B^\frac12(P^t-I)P^{KM-1-j}\|^2\E[\|A\Delta_j^\delta\|^2]\\
    &\leq c_0^2L\|B\|\|(P^t-I)P^{KM-1-j}\|^2\E[\|A\Delta_j^\delta\|^2],\\
    \E[T_{t,j}^2]& \leq c_0^2L \|B^\frac12 P^{KM+t-1-j}\|^2\E[\|A\Delta_j^\delta\|^2]\leq c_0^2 L\|B\| \E[\|A\Delta_j^\delta\|^2].
\end{align*}
Then, by Lemma \ref{lem:kernel}, we deduce  
\begin{align*}
{\rm I}_1
\leq & c_0^2 L\|B\| t^2\sum_{j=1}^{(K-1)M+t-2} (KM+t-1-j)^{-2} \E[\|A \Delta_j^\delta\|^2],\\
{\rm I}_2\leq & c_0^2 L\|B\| \sum_{j=(K-1)M+t-1}^{KM-1}\E[\|A \Delta_j^\delta\|^2],\\
{\rm I}_3\leq & c_0^2 L\|B\|\sum_{j=KM+1}^{KM+t-1} \E[\|A \Delta_j^\delta\|^2].
\end{align*}
Meanwhile, by the commutativity of $B,P$ and Lemma \ref{lem:kernel} with $\epsilon=\frac12$, we can alternatively get
\begin{align*}
{\rm I}_1
\leq &  4c_0 L t^2\sum_{j=1}^{(K-1)M+t-2} (KM+t-1-j)^{-3} \E[\|A \Delta_j^\delta\|^2],\\
{\rm I}_2
\leq &  \frac{c_0L}{2} \left(2c_0\|B\|\E[\|A \Delta_{KM-1}^\delta\|^2]+\sum_{j=(K-1)M+t-1}^{KM-2}(KM-1-j)^{-1}\E[\|A \Delta_j^\delta\|^2]\right),\\
{\rm I}_3
\leq &  \frac{c_0L}{2}\left(2c_0\|B\|\E[\|A \Delta_{KM+t-1}^\delta\|^2]+\sum_{j=KM+1}^{KM+t-2} (KM+t-1-j)^{-1}\E[\|A\Delta_j^\delta\|^2]\right).
\end{align*}
Similarly, when $K=0$, there hold 
\begin{align*}
&\E[\|A\Delta_{t}^\delta-\E[A\Delta_{t}^\delta]\|^2]
= \sum_{j=1}^{t-1}\E[ T_{t,j}^2]\\
\leq& \left\{
\begin{aligned}c_0^2 L\|B\|\sum_{j=1}^{t-1} \E[\|A \Delta_j^\delta\|^2],\\
\frac{c_0L}{2}\left(2c_0\|B\|\E[\|A \Delta_{t-1}^\delta\|^2]+\sum_{j=1}^{t-2} (t-1-j)^{-1} \E[\|A \Delta_j^\delta\|^2]\right).
\end{aligned}\right.
\end{align*}
Then combining the preceding estimates with the identity $\|B\|=n^{-1}\|A\|^2$ and the inequality $t\leq M-1$ gives the estimate \eqref{eqn:Delta_E2}.
Finally, when $K\geq1$, by Lemma \ref{lem:bias-var} and the triangle inequality, we derive
\begin{align*}
\|A\Delta_{KM+t}^\delta-\E[A\Delta_{KM+t}^\delta]\|
\leq&\sum_{j=1}^{KM-1} S_{t,j}+\sum_{j=KM+1}^{KM+t-1}T_{t,j}.
\end{align*}
Thus for any $\epsilon\in(0,\frac12]$, by Lemmas \ref{lem:kernel} and \ref{lem:N} and the identity $\|B\|=n^{-1}\|A\|^2$, we have
\begin{align*}
S_{t,j}
\leq & c_0\sqrt{nL} \|B^\frac12 (P^t-I)P^{KM-1-j}\|\|A\Delta_j^\delta\|\\
\leq&   \left\{
\begin{aligned}2^{1+\epsilon}\epsilon^\epsilon c_0^{1-\epsilon}n^{\epsilon}\sqrt{L}\|A\|^{1-2\epsilon}t(KM+t-1-j)^{-(1+\epsilon)}\|A\Delta_j^\delta\|, \quad \forall j\leq (K-1)M+t-2\\
c_0\sqrt{L} \|A\|\|A\Delta_j^\delta\|, \quad \forall j\geq (K-1)M+t-1\\
\end{aligned}\right.,\\
T_{t,j}\leq &c_0\sqrt{L}\|A P^{KM+t-1-j} \|\|A\Delta_j^\delta\|
\leq c_0\sqrt{L}\|A\|\|A\Delta_j^\delta\|.
\end{align*}
When $K=0$, there holds
\begin{align*}
&\|A\Delta_{t}^\delta-\E[A\Delta_{t}^\delta]\|
\leq \sum_{j=1}^{t-1}T_{t,j}
\leq c_0\sqrt{L}\|A\| \sum_{j=1}^{t-1} \|A \Delta_j^\delta\|.
\end{align*}
Combining these estimates with the inequality $t\leq M-1$ gives the estimate \eqref{eqn:Delta_E3}.
\end{proof}

The next lemma gives several basic estimates on the following summations:
\begin{align*}
\overline{\Phi}_{j_1}^{j_2}(0,1)&=2c_0\|B\|\E[\|A \Delta_{j_2}^\delta\|^2]+\sum_{j=j_1}^{j_2-1}(j_2-j)^{-1}\E[\|A \Delta_j^\delta\|^2],\\
\Phi_{j_1}^{j_2}(0,\tfrac12)&=\sqrt{2c_0\|B\|}\|A \Delta_{j_2}^\delta\|+\sum_{j=j_1}^{j_2-1}(j_2-j)^{-\frac12}\|A \Delta_j^\delta\|,\\
\overline{\Phi}_{j_1}^{j_2}(i,r)&=\sum_{j=j_1}^{j_2}(j_2+i-j)^{-r}\E[\|A \Delta_j^\delta\|^2],\\ \Phi_{j_1}^{j_2}(i,r)&=\sum_{j=j_1}^{j_2}(j_2+i-j)^{-r}\|A \Delta_j^\delta\|,
\end{align*}
for any $r,\;i>0$ or $r=0$. The proof of these results is elementary but lengthy and tedious in order to get precise constants.

\begin{lemma}\label{lem:sum}
For any $k\geq 1$, let $k_{\rm c}:=k-M-2$ and $k_M:=[\frac{k}{M}]M$. If there  holds
\begin{equation} \label{eqn:bdd-ADelta}
\max\big(\E[\|A\Delta_{k}^\delta\|^2]^\frac12,\|A\Delta_{k}^\delta\|\big)\leq c(k+M)^{-1},
\end{equation}  
then for any $k> M$, $k\neq k_M$ and $\epsilon\in(0,\frac12]$, there hold
\begin{align}
&\overline{\Phi}_{k_{\rm c}+1}^{k-1}(0,0)\leq c_{K,1}^2c^2 M (k+M)^{-2},\quad 
\Phi_{k_{\rm c}+1}^{k-1}(0,0)\leq c_{K,1} c M (k+M)^{-1},\label{eqn:sum00}\\
&\overline{\Phi}_{k_{\rm c}+1}^{k_M-1}(0,1)+
\overline{\Phi}_{k_M+1}^{k-1}(0,1)\leq  c_{K,1}^2 c^2 (1+2\ln M+4c_0\|B\|)(k+M)^{-2} ,\label{eqn:sum01}\\
&\overline{\Phi}_{1}^{k_{\rm c}}(M+1,2)\leq 4c_{K,1}^{2}c^2 M^{-1}(k+M)^{-2}, \label{eqn:sumM2}\\
&\overline{\Phi}_{1}^{k_{\rm c}}(M+1,3)\leq c_{K,1}^{2}c_{K,2}c^2 M^{-2}(k+M)^{-2}, \label{eqn:sumM3}\\
&\Phi_{1}^{k_{\rm c}}(M+1,\epsilon+1)\leq 2c_{K,1}c_{K,\epsilon}c M^{-\epsilon}(k+M)^{-1}, \label{eqn:sumMe}\\
&\overline{\Phi}_{1}^{k-1}(0,1)\leq (3+2c_0\|B\|)c^2 k^{-1}, \label{eqn:sum01_e}\\
&\Phi_{1}^{k-1}(0,\tfrac12)\leq 3\sqrt{2}c k^{-\frac12}\ln k, \label{eqn:sum00.5_e}
\end{align}
with $K=[\frac{k}{M}]$, and the constants $c_{K,1}$, $c_{K,2}$ and $c_{K,\epsilon}$ defined respectively by 
\begin{align*} 
c_{K,1}&=1+\frac{2}{K}\leq 3,\\ c_{K,2}&=2+\frac{3}{K+1}+\frac{6\ln M}{(K+1)^2}\leq \frac72+\frac32\ln M,\\ c_{K,\epsilon}&=\frac{e^{-1}+1}{\epsilon}+\frac{2^\epsilon \ln M}{(K+1)^{\epsilon}}.
\end{align*} 
Note that there hold $\lim_{K\to\infty}c_{K,1}=1$, $\lim_{K\to\infty}c_{K,2}=2$, and $\lim_{K\to\infty}c_{K,\epsilon}=(e^{-1}+1)\epsilon^{-1}$. 
\end{lemma}
\begin{proof}
Let $k=KM+t$ with $K\geq 1$ and $t=1,\cdots,M-1$. In the proof, we use the shorthand notation $\overline{k}=k-1$ and $\overline{k_M}=k_M-1$. Then there holds the elementary inequality:
\begin{equation}\label{eqn:kM}
    \overline{k}^{-1} \leq c_{K,1}(k+M)^{-1}.
\end{equation}
The estimates in \eqref{eqn:sum00} follow directly from  \eqref{eqn:bdd-ADelta}, the identity $\Delta_{KM}^\delta=\Delta_{(K-1)M}^\delta=0$ and \eqref{eqn:kM}:
\begin{align*}
\overline{\Phi}_{k_{\rm c}+1}^{\overline{k}}(0,0)&=\sum_{j=k_{\rm c}+1}^{\overline{k}}\E[\|A \Delta_j^\delta\|^2] 
\leq c^2 M\overline{k}^{-2}
\leq c_{K,1}^{2} c^2 M (k+M)^{-2},\\
\Phi_{k_{\rm c}+1}^{\overline{k}}(0,0)&=\sum_{j=k_{\rm c}+1}^{\overline{k}}\|A \Delta_j^\delta\|\leq c M \overline{k}^{-1}\leq c_{K,1}c M (k+M)^{-1}.
\end{align*}
Next for the estimate \eqref{eqn:sum01}, we have
\begin{align*}
&\overline{\Phi}_{k_{\rm c}+1}^{\overline{k_M}}(0,1)+
\overline{\Phi}_{k_M+1}^{\overline{k}}(0,1)\\=&\sum_{j=k_{\rm c}+1}^{k_M-2}(\overline{k_M}-j)^{-1}\E[\|A \Delta_j^\delta\|^2]+\sum_{j=k_M+1}^{k-2}(\overline{k}-j)^{-1}\E[\|A \Delta_j^\delta\|^2]\\
&+2c_0\|B\|\big(\E[\|A \Delta_{\overline{k_M}}^\delta\|^2]+\E[\|A \Delta_{\overline{k}}^\delta\|^2]\big)\\
=& \sum_{j=1}^{M-t}j^{-1}\E[\|A \Delta_{\overline{k_M}-j}^\delta\|^2]
+ \sum_{j=1}^{t-2} j^{-1}\E[\|A\Delta_{\overline{k}-j}^\delta\|^2]\\
&+2c_0\|B\|\big(\E[\|A \Delta_{\overline{k_M}}^\delta\|^2]+\E[\|A \Delta_{\overline{k}}^\delta\|^2]\big)\\
\leq&c^2\bigg(\sum_{j=1}^{M-t}j^{-1}(\overline{k_M}-j+M)^{-2}
+\sum_{j=1}^{t-2} j^{-1}(\overline{k}-j+M)^{-2}\\
&\quad+2c_0\|B\|\big((\overline{k_M}+M)^{-2}+(\overline{k}+M)^{-2}\big)\bigg)\\
\leq&c^2\overline{k}^{-2}{\rm I},
\end{align*} 
with 
\begin{align*} 
{\rm I}=  4c_0\|B\|+\sum_{j=1}^{M-t}j^{-1}
+\sum_{j=1}^{t-2} j^{-1}.
\end{align*}
Note that the factor ${\rm I}$ can be bounded by
\begin{align*}
{\rm I}&\leq  4c_0\|B\|+2+\ln(M-t)
+\ln t\leq 4c_0\|B\|+2+2\ln\tfrac{M}{2} \\
&= 4c_0\|B\|+2- 2\ln 2+2\ln M\leq 1+2\ln M+4c_0\|B\|.
\end{align*}
Then with the estimate \eqref{eqn:kM}, there holds
\begin{align*}
\overline{\Phi}_{k_{\rm c}+1}^{\overline{k_M}}(0,1)+
\overline{\Phi}_{k_M+1}^{\overline{k}}(0,1)\leq  c_{K,1}^2 c^2 (1+2\ln M+4c_0\|B\|)(k+M)^{-2}.
\end{align*}
Next, we derive the estimates \eqref{eqn:sumM2},  \eqref{eqn:sumM3} and \eqref{eqn:sum01_e}. 
For the estimate \eqref{eqn:sumM2}, by the splitting 
$$(j'-j)^{-2}j^{-2}=(j')^{-2}\big((j'-j)^{-1}+j^{-1}\big)^2\leq 2(j')^{-2}\big((j'-j)^{-2}+j^{-2}\big),$$ we obtain 
\begin{align*}
&\overline{\Phi}_{1}^{k_{\rm c}}(M+1,2) \leq c^2\sum_{j=1}^{k-M-2} (\overline{k}-j)^{-2} (j+M)^{-2}
= c^2\sum_{j=M+1}^{k-2} (\overline{k}+M-j)^{-2} j^{-2}\\
\leq & 2c^2(\overline{k}+M)^{-2}\sum_{j=M+1}^{k-2} \big((\overline{k}+M-j)^{-2}+ j^{-2}\big)
\leq 4c_{K,1}^{2}c^2M^{-1}(k+M)^{-2}.
\end{align*}
Likewise, for the estimate \eqref{eqn:sumM3}, by the splitting 
\begin{align*}
(j'-j)^{-3}j^{-2}=&3(j')^{-4}\big((j'-j)^{-1}+j^{-1}\big)\\
&+(j')^{-3}\big(2(j'-j)^{-2}+j^{-2}\big)+(j')^{-2}(j'-j)^{-3},
\end{align*}
we derive
\begin{align*}
&\overline{\Phi}_{1}^{k_{\rm c}}(M+1,3) \leq c^2\sum_{j=1}^{k-M-2} (\overline{k}-j)^{-3} (j+M)^{-2}
= c^2\sum_{j=M+1}^{k-2} (\overline{k}+M-j)^{-3} j^{-2}\\
=& c^2 \big(\overline{k}+M\big)^{-2}\bigg[6\big(\overline{k}+M\big)^{-2}\sum_{j=M+1}^{k-2}j^{-1}
+3\big(\overline{k}+M\big)^{-1}\sum_{j=M+1}^{k-2}j^{-2}+\sum_{j=M+1}^{k-2} j^{-3}\bigg]\\[2mm]
\leq& c^2\big(\overline{k}+M\big)^{-2}\Big[6\big((K+1)M\big)^{-2}\ln(KM+t-1)
+3\big((K+1)M\big)^{-1}M^{-1}+\tfrac12 M^{-2}\Big]\\[2mm]
\leq& \frac{c_{K,1}^{2}c^2}{(k+M)^{2}M^{2}}\Big[6(K+1)^{-2}\big(\ln(K+1)+\ln M\big) +\tfrac{3}{K+1}+\tfrac12 \Big].
\end{align*}
Then the inequality 
\begin{equation}\label{eqn:ln}
s^{-r}\ln s\leq (er)^{-1}, \quad \forall s,r > 0
\end{equation}
implies the desired bound on
$\overline{\Phi}_{1}^{k_{\rm c}}(M+1,3)$.
For the estimate \eqref{eqn:sum01_e}, the splitting
$$(j'-j)^{-1}j^{-2}=(j')^{-1}j^{-2}+(j')^{-2}\big((j'-j)^{-1}+j^{-1}\big)$$ implies
\begin{align*}
&\overline{\Phi}_{1}^{\overline{k}}(0,1) \leq  c^2\Big(2c_0\|B\| (\overline{k}+M)^{-2}+\sum_{j=1}^{k-2} (\overline{k}-j)^{-1}(j+M)^{-2}\Big)\\
=& c^2\Big(2c_0\|B\| (\overline{k}+M)^{-2}+\sum_{j=M+1}^{k+M-2} (\overline{k}+M-j)^{-1}j^{-2}\Big)\\
= &c^2(\overline{k}+M)^{-1}\bigg[\sum_{j=M+1}^{k+M-2}j^{-2}
+ (\overline{k}+M)^{-1}\Big(2c_0\|B\|
+\sum_{j=M+1}^{k+M-2} \big((\overline{k}+M-j)^{-1}+j^{-1}\big)\Big)\bigg]\\[2mm]
\leq&c^2(\overline{k}+M)^{-1}\big[M^{-1}+(\overline{k}+M)^{-1}\big(1+2\ln(\overline{k}+M)+2c_0\|B\|\big) \big].
\end{align*}
Then, using the inequality \eqref{eqn:ln}, we derive 
\begin{align*}
\overline{\Phi}_{1}^{\overline{k}}(0,1)
\leq c^2(\overline{k}+M)^{-1}\big((2+2c_0\|B\|)M^{-1}+2e^{-1}\big)
\leq (3+2c_0\|B\|)c^2 k^{-1}.
\end{align*}
Now, we derive the estimates \eqref{eqn:sum00.5_e} and \eqref{eqn:sumMe} by splitting the summations into two parts. Let $\widetilde{k_M} = (k+M-1)/2$. 
For the estimate \eqref{eqn:sum00.5_e}, there holds
\begin{align*}
\Phi_{1}^{\overline{k}}(0,\tfrac12) \leq  &\;c\Big(\sqrt{2c_0\|B\|}(\overline{k}+M)^{-1}+\sum_{j=1}^{k-2} (\overline{k}-j)^{-\frac12}(j+M)^{-1}\Big)\\
=&c\Big(\sqrt{2c_0\|B\|}(\overline{k}+M)^{-1}+\sum_{j=M+1}^{k+M-2} (\overline{k}+M-j)^{-\frac12}j^{-1}\Big)\\
\leq& c({\rm I}_{11}+{\rm I}_{12}),\\
\end{align*}
with 
\begin{align*}
{\rm I}_{11}=\sum_{j=M+1}^{[\widetilde{k_M}]} \widetilde{k_M}^{-\frac12}j^{-1}
\quad\mbox{and}\quad
{\rm I}_{12}=\bigg(\sqrt{2^{-1}c_0\|B\|}+\sum_{j=[\widetilde{k_M}]+1}^{2\widetilde{k_M}-1} (\overline{k}+M-j)^{-\frac12}\bigg)\widetilde{k_M}^{-1}.
\end{align*}
The decomposition is well-defined with the convention $\sum_{j=i}^{i'}R_j=0$ for any $\{R_j\}_j$ and $i'< i$. Then we have
\begin{equation*} 
{\rm I}_{11}
\leq \widetilde{k_M}^{-\frac12}\ln \widetilde{k_M}
\leq\sqrt{2}k^{-\frac12}\ln k\quad \mbox{and}\quad 
{\rm I}_{12}
\leq 2\widetilde{k_M}^{-\frac12}
\leq 2\sqrt{2}k^{-\frac12}.
\end{equation*}
Similarly, for the estimate \eqref{eqn:sumMe}, when $\epsilon\in(0,\frac12]$, we split $\Phi_{1}^{k_{\rm c}}(M+1,\epsilon+1)$ into  
\begin{align*}
&\Phi_{1}^{k_{\rm c}}(M+1,\epsilon+1)
\leq  c\sum_{j=M+1}^{k-2} \big(\overline{k}+M-j\big)^{-(1+\epsilon)}j^{-1}
\leq c({\rm I}_{21}+{\rm I}_{22}),
\end{align*}
with 
\begin{align*} 
{\rm I}_{21}=\sum_{j=M+1}^{[\widetilde{k_M}]} \widetilde{k_M}^{-(1+\epsilon)}j^{-1} \quad\mbox{and}\quad {\rm I}_{22}=\sum_{j=[\widetilde{k_M}]+1}^{k-2} \big(\overline{k}+M-j\big)^{-(1+\epsilon)}\widetilde{k_M}^{-1}.
\end{align*}
Then
\begin{align*}
{\rm I}_{21}&
\leq \widetilde{k_M}^{-(1+\epsilon)}\ln\widetilde {k_M}
\leq2\Big((e\epsilon)^{-1}+2^\epsilon(K+1)^{-\epsilon} \ln M\Big)M^{-\epsilon}(\overline{k}+M)^{-1},\\[1mm]
{\rm I}_{22}&
\leq 2\epsilon^{-1} M^{-\epsilon}(\overline{k}+M)^{-1}.
\end{align*}
Finally, the inequality $(\overline{k}+M)^{-1}\leq c_{K,1} (k+M)^{-1}$
completes the proof of the lemma.
\end{proof}

The proof uses also the following elementary estimate on the function 
\begin{equation}\label{eqn:f}
f(\epsilon)=2^{2+\epsilon}\epsilon^\epsilon c_0^{1-\epsilon}n^{\epsilon}\sqrt{L} M^{1-\epsilon}\|A\|^{1-2\epsilon}c_{K,\epsilon}.
\end{equation}
\begin{lemma}\label{lem:f}
If $c_0<C_0$ and $K\geq M-1$, then $\inf_{\epsilon\in (0,\frac12]}f(\epsilon)\leq \frac{4\sqrt{e}}{7}$.
\end{lemma}
\begin{proof}
By the definition of $f(\epsilon)$ and the inequality \eqref{eqn:ln}, we have
\begin{align*}
f(\epsilon)
=&2^{2+\epsilon}\epsilon^\epsilon c_0^{1-\epsilon}n^{\epsilon}\sqrt{L} M^{1-\epsilon}\|A\|^{1-2\epsilon}\big((e^{-1}+1)\epsilon^{-1}+2^\epsilon(K+1)^{-\epsilon} \ln M\big)\\
=&4\big(c_0\sqrt{L}M\|A\|\epsilon^{-1}(2n\sqrt{L}\|A\|^{-1})^\frac{\epsilon}{1-\epsilon}\big)^{1-\epsilon}\big(e^{-1}+1+2^\epsilon \epsilon (K+1)^{-\epsilon} \ln M\big)\\
\leq &8\big(c_0\sqrt{L}M\|A\|\epsilon^{-1}(2n\sqrt{L}\|A\|^{-1})^\frac{\epsilon}{1-\epsilon}\big )^{1-\epsilon},
\end{align*}
for any $\epsilon\in(0,\frac12]$ and $K\geq M-1$. Let $g(\epsilon)=\epsilon^{-1}(2n\sqrt{L}\|A\|^{-1})^\frac{\epsilon}{1-\epsilon}$. Then 
\begin{align*}
g'(\epsilon)=\epsilon^{-2}(1-\epsilon)^{-2}(2n\sqrt{L}\|A\|^{-1})^\frac{\epsilon}{1-\epsilon}\big(\epsilon\ln(2n\sqrt{L}\|A\|^{-1})-(1-\epsilon)^{2}\big).
\end{align*}
The fact $\|A\|\leq \sqrt{ \sum_{i=1}^n\|A_i\|^2}\leq \sqrt{nL}$ implies $$c:=2+\ln(2n\sqrt{L}\|A\|^{-1})\geq 2+\ln (2\sqrt{n})>2+\ln2.$$  $g(\epsilon)$ attains its minimum over the interval $(0,\frac12]$ at $\epsilon=\epsilon^*=\frac{2}{c+\sqrt{c^2-4}}<\frac12$, and $$g(\epsilon^*)=\frac{c+\sqrt{c^2-4}}{2}
e^{\frac{2(c-2)}{c+\sqrt{c^2-4}-2}}\leq c e.$$
Thus, for $c_0<C_0$, we have
\begin{align*}
f(\epsilon^*)\leq &8\big(c_0\sqrt{L}M\|A\|ce\big )^{1-\frac{2}{c+\sqrt{c^2-4}}}\\
\leq& 8\Big(ec_0\sqrt{L}M\|A\|\ln(2e^2n\sqrt{L}\|A\|^{-1})\Big )^{\frac12}<\tfrac{4\sqrt{e}}{7}.
\end{align*}
This completes the proof of the lemma.
\end{proof}

Now we can prove Theorem \ref{thm:Delta} by mathematical induction. 
\begin{proof}
If $k\leq K_0M$ with some $K_0\geq 1$, the estimate \eqref{eqn:Delta_thmE1} holds for any sufficiently large $c_1$ and $c_2$. Now assume that it
holds up to $k=KM+t-1$ with some $K\geq K_0$ and $1\leq t\leq M-1$. Then we prove the assertion for the case $k=KM+t$. (It holds trivially when $t=0$, since $\Delta_{KM}^\delta=0$.)
Fix $k=KM+t$ and let $k_{\rm c}:=k-M-2$ and $k_M:=[\frac{k}{M}]M=KM$.
By the bias-variance decomposition,
and the estimates \eqref{eqn:Delta_E1} and \eqref{eqn:Delta_E2}
in Lemma \ref{lem:decom_Delta}, we have
\begin{align*}
\E[\|A\Delta_{k}^\delta\|^2]\leq& 2(\|A\|^2\|e_0^\delta\|^2+ \delta^2) M^2 k^{-2}\\
&+n^{-1}c_0^2 L\|A\|^2\Big(M^2\overline{\Phi}_{1}^{k_{\rm c}}(M+1,2)+\overline{\Phi}_{k_{\rm c}+1}^{k-1}(0,0)\Big).
\end{align*}
Then, by setting $c=\sqrt{c_1+c_2\delta^2}$, the estimates \eqref{eqn:sum00} and \eqref{eqn:sumM2} and the inequality $k^{-1}\leq c_{K,1}(k+M)^{-1}$ (cf. \eqref{eqn:kM}) with $c_{K,i}$ given in Lemma \ref{lem:sum} yield
\begin{align*}
\E[\|A\Delta_{k}^\delta\|^2]\leq& c_{K,1}^{2}\Big[5c_0^2LMn^{-1}\|A\|^2(c_1+c_2\delta^2)\\
&+2(\|A\|^2\|e_0^\delta\|^2+ \delta^2) M^2\Big](k+M)^{-2}\\
\leq &(c_1+c_2\delta^2)(k+M)^{-2},
\end{align*}
for any $c_0<3^{-1}\|A\|^{-1}(5Ln^{-1}M)^{-\frac12}$ and $K\geq 1$, with sufficiently large $c_1$ and $c_2$. Alternatively, using the second estimate in \eqref{eqn:Delta_E2}, we can bound $\E[\|A\Delta_{k}^\delta\|^2]$ by
\begin{align*}
\E[\|A\Delta_{k}^\delta\|^2]\leq& 2(\|A\|^2\|e_0^\delta\|^2+ \delta^2) M^2 k^{-2}\\
&+\tfrac{c_0 L}{2}\Big(8M^2\overline{\Phi}_{1}^{k_{\rm c}}(M+1,3)+\overline{\Phi}_{k_{\rm c}+1}^{k_M-1}(0,1)+\overline{\Phi}_{k_M+1}^{k-1}(0,1)\Big).
\end{align*}
Then, using the estimates \eqref{eqn:sum01} and \eqref{eqn:sumM3}, we derive
\begin{align*}
\E[\|A\Delta_{k}^\delta\|^2]\leq& c_{K,1}^2\big[c_0L(\tfrac12+\ln M+2c_0\|B\|+4
c_{K,2})(c_1+c_2\delta^2)\\
&+2(\|A\|^2\|e_0^\delta\|^2+ \delta^2) M^2\big](k+M)^{-2}\\
\leq& (c_1+c_2\delta^2)(k+M)^{-2},
\end{align*}
for any $c_0<9^{-1}L^{-1}(15+7\ln M)^{-1}$ and $K\geq 1$, with sufficiently large $c_1$ and $c_2$.
This completes the proof of the estimate \eqref{eqn:Delta_thmE1}.

\noindent Next, we prove the estimate \eqref{eqn:Delta_thmE2}.
Similarly, for the cases $k\leq K_0M$ with some $K_0\geq 1$, the estimate \eqref{eqn:Delta_thmE2} holds trivially for sufficiently large $c_1$ and $c_2$. Now, assume that the bound
holds up to $k=KM+t-1$ with some $K\geq K_0$ and $1\leq t\leq M-1$, and prove the assertion for the case $k=KM+t$. Fix $k=KM+t$ and let $k_{\rm c}:=k-M-2$ and $k_M:=[\frac{k}{M}]M=KM$.
By the triangle inequality
$\|A\Delta_{k}^\delta\|\leq \|\E[A\Delta_{k}^\delta]\|+\|\Delta_{k}^\delta-\E[A\Delta_{k}^\delta]\|,$ and \eqref{eqn:Delta_E1} and \eqref{eqn:Delta_E3}, we have
\begin{align}
&\|A\Delta_{k}^\delta\|\leq (\|A\|\|e_0^\delta\|+ \delta) M k^{-1}+\Big({\rm I}_1+c_0\sqrt{L}\|A\|\Phi_{k_{\rm c}+1}^{k-1}(0,0)\Big),\label{eqn:Delta_as}
\end{align}
with
$$
{\rm I}_1=2^{1+\epsilon}\epsilon^\epsilon c_0^{1-\epsilon}n^{\epsilon}\sqrt{L}\|A\|^{1-2\epsilon}M \Phi_{1}^{k_{\rm c}}(M+1,1+\epsilon).
$$ 
By \eqref{eqn:sumMe} (with $c=c_1+c_2\delta$), we derive
\begin{align*}
{\rm I}_1\leq& 2^{2+\epsilon}\epsilon^\epsilon (c_1+c_2\delta)c_0^{1-\epsilon}n^{\epsilon}\sqrt{L}\|A\|^{1-2\epsilon} M^{1-\epsilon}c_{K,1}c_{K,\epsilon}(k+M)^{-1}\\
=&(c_1+c_2\delta)f(\epsilon) c_{K,1}(k+M)^{-1},
\end{align*}
with $f(\epsilon)$ given in \eqref{eqn:f}.
This, \eqref{eqn:Delta_as}, \eqref{eqn:sum00} in Lemma \ref{lem:sum}, and the inequality \eqref{eqn:kM} yield
\begin{equation}\label{eqn:Delta_as0}
\begin{array}{cc}
\|A\Delta_{k}^\delta\|\leq c_{K,1}\Big[(c_1+c_2\delta)\big(c_0\sqrt{L}M\|A\|+f(\epsilon)\big)+(\|A\|\|e_0^\delta\|+ \delta) M \Big](k+M)^{-1}.
\end{array}
\end{equation}
Then by Lemma \ref{lem:f} and the inequality $c_0\sqrt{L}M\|A\|<\big(14^2(2+\ln 2)\big)^{-1}$ when $c_0<C_0$, we derive from \eqref{eqn:Delta_as0} that
\begin{align*}
\|A\Delta_{k}^\delta\|\leq (c_1+c_2\delta)(k+M)^{-1},
\end{align*}
for any $K\geq 35$, with sufficiently large $c_1$ and $c_2$.
This completes the proof of the theorem.
\end{proof}

\bibliographystyle{abbrv}
\bibliography{sgd}
\end{document}